\documentclass[11pt]{amsart}
\usepackage{url}
\usepackage{hyperref}
\usepackage[mathscr]{eucal}
\usepackage{amssymb}

\newtheorem{theorem}{Theorem}[section]
\newtheorem{lemma}[theorem]{Lemma}
\newtheorem{proposition}[theorem]{Proposition}
\newtheorem{corollary}[theorem]{Corollary}

\newtheorem{remark}[theorem]{Remark}
\newtheorem{remarks}[theorem]{Remarks}
\def\notdiv{\nmid}

\def\lien{\mathrel{\mkern-4mu}}
\def\too{\relbar\lien\rightarrow}
\def\tooo{\relbar\lien\relbar\lien\too}

\def\Q{\mathbb{Q}}
\def\Z{\mathbb{Z}}
\def\F{\mathbb{F}}
\let\ds=\displaystyle
\let\ov=\overline
\let\wt=\widetilde
\let\es=\emptyset
\def\Cl{{\mathcal C}\hskip-2pt{\ell}}
\def\cl{c\hskip-1pt{\ell}}
\def\order{\raise1.5pt \hbox{${\scriptscriptstyle \#}$}}
\def\prd{ \ds\mathop{\raise 2.0pt \hbox{$\prod$}}\limits}
\def\sm{ \ds\mathop{\raise 2.0pt \hbox{$\sum$}}\limits}
\def\s{\hskip1pt}

\begin{document}

\markboth{Georges Gras}{}

\title[Note on some $p$-invariants of $\Q\big(N^\frac{1}{p}_{} \big)$]
{Note on some $p$-invariants of $\Q\big(N^\frac{1}{p}\big)$ \\
using reflection theorem \\
\footnotesize Tables of the groups $\Cl$, ${\mathcal T}$, $\wt \Cl$}


\author{Georges Gras}
\address{Villa la Gardette, 4 Chemin Ch\^ateau Gagni\`ere,
F-38520 Le Bourg d'Oisans}
\email{g.mn.gras@wanadoo.fr}

\date{September 30, 2019}

\begin{abstract} 
Let $p > 2$ be a prime number and let 
$N=p^\delta \prod_{i=1}^n \ell_i^{\lambda_i}$, with primes $\ell_i$,
$\delta \geq 0$, $\lambda_i \not\equiv 0 \pmod p$.
We consider the $p$-class groups $\Cl_L$, $\Cl_M$ of the fields
$L:=\Q(N^\frac{1}{p})$ and $M :=\Q(N^\frac{1}{p}, \mu_p)$, by comparison 
with the $p$-torsion groups ${\mathcal T}_L$ and ${\mathcal T}_M$
of the abelian $p$-ramification theory, in the framework 
of the reflection theorem, and obtain relations between the ranks 
of the isotypic components (Theorem \ref{thmrefl}). 

\noindent
For $p=3$, we characterize the integers $N$ such that 
$L=\Q(N^\frac{1}{3})$ is $3$-rational (i.e., ${\mathcal T}_L=1$), 
giving the following values: $N=3$; $N=3^\delta \ell$, $\ell = -1+ 3\, u$;
$N=3^\delta \ell$, $\ell = {\rm N}_{\Q( \sqrt{-3})/\Q}
\big(1+3\, \big (\frac{a+b \sqrt{-3}}{2} \big) \big)$, 
with $\ell$ prime and $3 \nmid u \, a \, b$ (Theorem \ref{3rat}).
We show that the $3$-class group $\Cl_L$ of $L$ is trivial if and only 
$\Cl_M=1$ (Theorem \ref{class}).
We give various tables with PARI/GP programs computing the 
structure of $\Cl_L$, $\Cl_M$, ${\mathcal T}_L$, ${\mathcal T}_M$
and of the logarithmic class groups (Appendix \ref{A}, \ref{B}, \ref{C}). 
\end{abstract}

\maketitle

\tableofcontents

\section{Introduction -- Literature} 

Let $p>2$ be prime. Many papers have studied the $p$-class groups 
$\Cl_L$ and $\Cl_M$ of non-Galois metabelian fields $L$ (like 
$L=\Q(\sqrt[p]{N})$ and its Galois closure $M =\Q({\sqrt[p]{N}, \mu_p})$) 
(e.g., \cite{Gerth, Gras0, Kobayashi} for $p=3$ then 
\cite{Iimura, Jaulent} in a general setting). For $p=3$
recent papers give extensive results on the $3$-structure of 
the class group $\Cl_M$ of $M$ \cite{AAIMT,ATIA1,ATIA2,HW}, and
considered in \cite{Schoof} with a new algebraic approach for any $p$.

\smallskip
Other recent papers give deep information about $L=\Q(\sqrt[p]{N})$ 
when $N \equiv 1 \pmod p$ is a prime number, using Galois 
cohomology, modular forms \cite{Calegari-Emerton,Schaefer-Stubley}; 
in \cite{Lecouturier} these results are obtained from class field 
theory, $N$-adic Gamma function and Gauss sums. A typical result is 
that the $p$-rank of the class group of $L$ is $\geq 2$ if and only if
$C:=\prod_{k=1}^{(N-1)/2} k^k$ is a $p$th power modulo $N$
(equivalent to $\sum_{k=1}^{(N-1)/2} k \cdot {\rm lg}_p(k) \equiv 0
\pmod p$, where ${\rm lg}_p$ is the discrete logarithm relative
to the choice of a primitive root modulo $N$).

\smallskip
In the classical framework, the more precise results for an arbitrary integer 
$N$ are that of Iimura (e.g., \cite[Lemma 1.1, Corollary to Theorem 2.3]{Iimura}), 
after the general study by Jaulent \cite{Jaulent} of the arbitrary metabelian 
case introducing the class field theory material establishing the fundamental
relations between the isotypic components of the class groups.

\smallskip
In this paper we introduce the supplementary $p$-groups 
${\mathcal T}_L$, ${\mathcal T}_M$, $\wt \Cl_L$, $\wt \Cl_M$,
stemming from abelian $p$-ramification theory.

\medskip
{\bf Main results.}\,{\it We consider any integer 
$N=p^{\delta} \ell_1^{\lambda_1} \cdots \ell_n^{\lambda_n}$, 
$\delta \geq 0$, $\ell_i$ primes, $\lambda_i \not \equiv 0 \pmod p$,
and introduce the finite $p$-groups ${\mathcal T}_L$ and ${\mathcal T}_M$ 
which are the torsion groups of the Galois groups of the maximal abelian 
$p$-ramified pro-$p$-extensions of $L:=\Q(\sqrt[p]{N})$ and 
$M:=\Q(\sqrt[p]{N}, \mu_p)$, respectively, which gives relations 
with the $p$-class groups $\Cl_L$ and $\Cl_M$ from reflection theorem. 

\smallskip
(i) A typical result is the following relation (with obvious notations)
between the ranks of the isotopic components of these $p$-groups 
regarding the group ${\rm Gal}(M/L) \simeq (\Z/p\Z)^\times$ 
(Theorem \ref{thmrefl}):
$$ {\rm rk}_{\chi^*}({\mathcal T}_M)=
{\rm rk}_\chi(\Cl_M/\cl_M(P_M)) + d_p, $$

where $\chi^*=\omega\cdot \chi^{-1}$, $\omega$ being the 
Teichm\"uller character, $d_p=1$ (resp. $0$) if 
$p\nmid N$ and $N^{p-1} \equiv 1 \pmod {p^2}$ (resp. if not),
$P_M$ being the set of $p$-places of $M$ ($\order P_M=1$ if 
$d_p=0$, $\order P_M=p$ if $d_p=1$); which gives for instance:
$${\rm rk}_p({\mathcal T}_L)  = 
{\rm rk}_\omega(\Cl_M/\cl_M(P_M))  + d_p. $$
 
(ii) We recall some criteria of $p$-rationality of $M$
(i.e., ${\mathcal T}_M=1$) and of $p$-principality in Theorem \ref{equiv}.
As an effective result, we obtain in Theorem \ref{3rat} that $L$ is $3$-rational 
if and only if $M$ is $3$-rational, which holds if and only if $N$ is of the 
following forms, up to $\Q^{\times 3}$ (where $\ell$ is a prime number):

\smallskip
\quad $\bullet$ $N=3$;

\smallskip
\quad $\bullet$ $N=\ell$, $3\,\ell$, $9\,\ell$, with 
$\ell = -1+ 3\, u$, $u \not\equiv 0 \!\!\pmod 3$;

\smallskip
\quad $\bullet$ $N=\ell$, $3\,\ell$, $9\,\ell$, with $\ell =
 {\rm N}_{\Q(\sqrt{-3})/\Q} \big (1+3\, \big( \frac{a+b \sqrt{-3}}{2} \big) \big)$, 
$a, b \not\equiv 0 \!\!\pmod 3$.

\smallskip
Similarly, we prove in Theorem \ref{class} that the $3$-class group of $L$ 
is trivial if and only if that of $M$ is trivial,
which is equivalent to $n = \nu_3 + d_3$ where $\nu_3=0$
(resp. $\nu_3=1$) if $\ell_i \equiv \pm1 \pmod 9$ for all $i$ (resp. if not).

\smallskip
(iii) In this context, we take into account \S\,\ref{bp} the 
logarithmic class group $\wt \Cl_F$ of Jaulent \cite{Jaulent0, Jalog},
related to the previous groups by means of a surjective 
map $\wt \Cl_F \to \Cl_F/\cl_F(P_F)$, for any field $F$;
when $F$ contains $\mu_p$ a reflection theorem does exist 
between $\wt \Cl_F$ and the Bertrandias-Payan module 
${\mathcal T}_F^{\rm bp}$.}

\section{Class groups and $p$-torsion groups -- Reflection theorem}

\subsection{Diagram of $p$-ramification}\label{shema}
Consider an arbitrary number field $F$ and the following diagram under 
the Leopoldt conjecture for $p$ in $F$. We denote by $P_F$ the set of 
$p$-places of $F$ and by $H_F^{\rm pr}$ (resp. $H_F$) 
the maximal abelian $P_F$-ramified pro-$p$-extension
(resp. the $p$-Hilbert class field) of $F$.

\smallskip
Let $U_F :=\bigoplus_{{\mathfrak p}\, \mid \, p} U_{\mathfrak p}$, be the
product of the groups of principal local units of the completion 
$F_{\mathfrak p}$ of $F$, for ${\mathfrak p} \mid p$, and let $\ov E_F$ 
be the closure of the diagonal image, of the group of principal 
global units $E_F$, in $U_F$.

\smallskip
We denote by $W_F := \bigoplus_{{\mathfrak p}\, \mid\, p}\mu^{}_{F_{\mathfrak p}}$ 
the torsion group of the $\Z_p$-module $U_F$.

\smallskip
The following $p$-adic result is valid under Leopoldt's conjecture for $p$ in~$F$:

\begin{lemma} \label{exact}
We have the exact sequence (${\rm log}_p$ being the $p$-adic logarithm):
\begin{equation*}
\begin{aligned}
1\to {\mathcal W}_F := W_F \big /\mu_p(F) & \tooo
{\rm tor}_{\Z_p}^{} \big(U_F \big /\ov E_F \big) \\
&\mathop {\tooo}^{{\rm log}_{p}}  {\rm tor}_{\Z_p}^{}\big({\rm log}_p
\big (U_F \big) \big / {\rm log}_p (\ov E_F) \big) \to 0. 
\end{aligned}
\end{equation*}
and the following schema:
\unitlength=0.95cm 
$$\vbox{\hbox{\hspace{-2.5cm} 
 \begin{picture}(11.5,5.4)
\put(6.5,4.50){\line(1,0){1.3}}
\put(8.7,4.50){\line(1,0){2.1}}
\put(3.85,4.50){\line(1,0){1.4}}
\put(9.2,4.2){\footnotesize$\simeq\! {\mathcal W}_F$}
\put(4.2,2.50){\line(1,0){1.25}}
\bezier{350}(3.8,4.9)(7.6,5.6)(11.0,4.9)
\put(7.2,5.4){\footnotesize${\mathcal T}_F$}

\bezier{350}(3.7,4.2)(5.85,3.4)(8.0,4.2)
\put(6.0,3.5){\footnotesize${\mathcal T}_F^{\rm bp}$}
\put(3.50,2.9){\line(0,1){1.25}}
\put(3.50,0.9){\line(0,1){1.25}}
\put(5.7,2.9){\line(0,1){1.25}}
\bezier{300}(3.9,0.5)(4.7,0.5)(5.6,2.3)
\put(5.2,1.3){\footnotesize$\simeq \! \Cl_F$} 
\bezier{300}(6.3,2.5)(8.5,2.6)(10.8,4.3)
\put(8.4,2.7){\footnotesize$\simeq \! U_F/\ov {E\,}_{\!F}$}
\put(10.85,4.4){$H_F^{\rm pr}$}
\put(5.4,4.4){$\wt {F\,} \!H_F$}
\put(7.85,4.4){$H_F^{\rm bp}$}
\put(6.7,4.2){\footnotesize$\simeq\! {\mathcal R}_F$}
\put(4.25,4.15){\footnotesize$\simeq\! \ov \Cl_F$}
\put(3.3,4.4){$\wt {F\,}$}
\put(5.5,2.4){$H_F$}
\put(2.9,2.4){$\wt {F\,}\! \!\cap \! H_F$}
\put(3.4,0.40){$F$}
\put(8.9,1.5){\footnotesize ${\mathcal A}_F$}
\bezier{500}(3.9,0.4)(9.5,0.8)(11.0,4.3)
\end{picture}   }}$$
\end{lemma}

\unitlength=1.0cm
By definition, ${\mathcal T}_F := {\rm tor}_{\Z_p} 
\big ({\mathcal A}_F \big)$ is the Galois group 
${\rm Gal}(H_F^{\rm pr} / \wt {F\,})$, where $ \wt {F\,}$
is the compositum of the $\Z_p$-extensions of $F$; let
$\ov \Cl_F$ be the subgroup of $\Cl_F$ corresponding to
${\rm Gal}(H_F/ \wt {F\,} \! \cap H_F)$ by class field theory
and put ${\mathcal R}_F := {\rm Gal}(H_F^{\rm bp} /\wt {F\,} H_F)$.
Then from the schema we get:
\begin{equation}\label{partial}
\order {\mathcal T}_F=\big [H_F\! :\! \wt {F\,}\! \cap H_F \big]
 \cdot \order {\rm tor}_{\Z_p} \big( U_F \big / \ov {E\,}_{\!F} \big) 
=\order \ov \Cl_F
\cdot \order {\mathcal R}_F \cdot \order {\mathcal W}_F .
\end{equation}
Of course, for $p \geq p_0$ (explicit), $\order {\mathcal W}_F=
\ov \Cl_F =1$, whence ${\mathcal T}_F={\mathcal R}_F$.

\begin{remarks}\label{w}
(i) We have ${\rm Gal}(H_F^{\rm pr} / H_F) \simeq U_F/\ov {E\,}_{\!F}$,
in which the image of ${\mathcal W}_F$ fixes $H_F^{\rm bp}$,
the Bertrandias--Payan field, 
${\mathcal T}_F^{\rm bp} := {\rm Gal}(H_F^{\rm bp} / \wt {F\,})
\simeq {\mathcal T}_F/{\mathcal W}_F$ being 
the Bertrandias--Payan module on the $p$-cyclic 
embedding problem \cite{BP}. 
Then ${\mathcal R}_F \simeq {\rm Gal}(H_F^{\rm bp} / \wt {F\,}\! H_F)$,
the ``normalized regulator'' ${\mathcal R}_F$, is closely related to the 
classical $p$-adic regulator of $F$ (see \cite[Proposition 5.2]{Gras2}).

\smallskip
(ii) In the case of the fields $M=\Q(\sqrt[p]{N}, \mu_p)$ and 
$L=\Q(\sqrt[p]{N})$ for $p>2$, we have ${\mathcal W}_M \simeq (\Z/p \Z)^{p-1}$ 
(resp. $\{0\}$) if $p$ splits in $M/K$ (resp. if not); then ${\mathcal W}_{L}
\simeq \Z/p \Z$ (resp. $\{0\}$). 
For instance, in the following example for $p=7$, the $7$-rank of 
${\mathcal W}_M$ is $6$ since $7$ splits in $M/K$, which explains
the large rank of ${\mathcal T}_M$, but ${\mathcal W}_M$ is not 
necessarily a direct factor:

\smallskip
\scriptsize
\begin{verbatim}
N = 197 = [197, 1]
Structure of C_L = [7]          Structure of C_M = [7, 7, 7, 7, 7]
Structure of T_L = [49]         Structure of T_M = [49, 7, 7, 7, 7, 7, 7, 7, 7]
\end{verbatim}
\normalsize
\end{remarks}

\subsection{Explicit genera theory for the field $L$}
It is easy to find a lower bound for the $p$-rank of $\Cl_L$ using
explicit unramified degree $p$ cyclic extensions of $L$ 
(i.e., genera theory).

\begin{proposition}\label{r0}
Let $N = p^{\delta} \prod_{i=1}^n \ell_i^{\lambda_i}$,
$\delta \geq 0$, $n\geq 0$, $\lambda_i \not\equiv 0 \pmod p$.
We have ${\rm rk}_p(\Cl_L) \geq \order 
\big\{i \in [1,n], \  \ell_i \equiv 1 \pmod p \big\}.$
\end{proposition}

\begin{proof}
Let $\ell \equiv 1 \pmod p$ be such a prime divisor of $N$ and let
$F_\ell$ be the subfield of $\Q(\mu_\ell)$ of degree $p$.
The abelian extensions $M F_{\ell}/M$ and
$L F_{\ell}/L$ are unramified since $F_{\ell}/\Q$ is unramified at $p$
and since $\ell$ is tamely ramified in the extensions 
$M F_{\ell}/K$ and $L F_{\ell}/\Q$, giving the result since the
$F_{\ell_i}$, $\ell_i \equiv 1 \pmod p$, are linearly independent.
\end{proof}

Of course, the injective map $\Cl_L \to \Cl_M$ implies the same 
lower bound for $\Cl_M$, but the Chevalley formula in $M/K$ 
shows that $\Cl_M$ is non-trivial as soon as the number $t$ of 
prime ideals of $K$, ramified in $M/K$, fulfills the condition 
$t >\frac{p-1}{2}+1$.
More precisely $t=r_p + \sum_{\ell \mid N} \, n_\ell$, where $n_\ell$
is the number of prime ideals of $K$ over $\ell$ and $r_p=1$ 
(resp. $0$) if $p$ is ramified (resp. unramified) in $M/K$.

\subsection{Characters and isotypic components of $\Cl_M$}

We consider the automorphism group $g={\rm Gal}(M/L) \simeq 
{\rm Gal}(K/\Q)$ with its characters $\chi=\omega^i$,
$i \in [1,p-1]$.
We denote by $e_\chi$ the corresponding idempotents and by
${\rm rk}_\chi(A_M)$ the $\F_p$-dimension of $(A_M/A_M^p)^{e_\chi}$ for
any finitely generated $\Z_p[g]$-module $A_M$ attached to $M$. 

\smallskip
Note that for the unit character
$\chi_0^{}$, ${\rm rk}_{\chi_0^{}}(A_M)=
{\rm dim}_{\F_p}((A_M/A_M^p)^{e_{\chi_0^{}}})$ is the
$p$-rank ${\rm rk}_p(A_L)$ of $A_L$ since the norm in $M/L$ is surjective 
when one considers $A_M=\Cl_M, {\mathcal T}_M$, and so on.

\subsubsection{Jaulent's results}
We apply his results in the metabelian framework defined by $M/\Q$.
Let $G :={\rm Gal}(M/K) =: \langle \sigma \rangle \simeq \Z/p\Z$.
We have of course the relation $\Cl_M=\bigoplus_{\chi} \Cl_M^{e_\chi}$
and we put $\Cl_M^{e_\chi}=: \Cl_{M,\chi}$ or $\Cl_M^\chi$ if more 
convenient. 

\smallskip
Because of the relations $s \cdot \sigma \cdot s^{-1} = \sigma^{\omega(s)}$, 
for all $s \in g$, we have, with:
$$\Theta := \hbox{$\frac{1}{p-1}$} \sm_{s \in g} \omega (s^{-1}) \, \sigma^{\omega(s)},$$

the classical rule:
\begin{equation}\label{permut}
\Theta \cdot e_{\chi} = e_{\chi \omega} \cdot \Theta,\ \,
\hbox{for any $\chi$}\ \ \hbox{\cite[Proposition 2]{Jaulent}.}
\end{equation}

The main property
of the ${\rm Gal}(M/\Q)$-modules deduced from $\Cl_M$ is the 
exact sequence of $\Z_p[G]$-modules \cite[Theorem 1]{Jaulent}:
\begin{equation}\label{jfj}
1 \to \Cl_{M,\chi}^G \too \Cl_{M,\chi} 
\mathop{\tooo}^{\Theta} \Cl_{M,\chi \omega} 
\too {\mathcal G}_{M,\chi \omega} \to 1,
\end{equation}

where $\Cl_M^G$ is the group of ambiguous $p$-classes
(since $\Theta \Z_p[G] = (1-\sigma)\Z_p[G]$), and where
${\mathcal G}_M :=\Cl_M/\Cl_M^\Theta$ is the $p$-genus group
with $\order {\mathcal G}_{M} = \order \Cl_M^G$.

\smallskip
When $p$ is irregular, the $p$-groups $\Cl_{K,\chi}$
intervenne and the above exact sequence
leads to the class number formula \cite[Theorem 4]{Jaulent}.

\subsubsection{Iimura's results}
In this metabelian context, the classical filtration 
of $\Cl_M$ in the $p$-cyclic extension $M/K$ is the following.
For $0\leq i \leq p-1$, let:
$$\hbox{${\mathcal I}_{M,i} := \Cl_M^{\Theta^i} \Cl_M^p /
\Cl_M^{\Theta^{i+1}} \Cl_M^p$, \   and 
$r_i(L) :={\rm dim}_{\F_p}(({\mathcal I}_{M,i})^{e_{\chi_0}})$;} $$

then \cite[Lemma 1.1]{Iimura}:
\begin{equation}\label{sigma}
{\rm rk}_p(\Cl_L)= \sm_{i=0}^{p-2}\  r_i(L).
\end{equation}

Moreover  \cite[Lemma 1.2]{Iimura}, we know that
$r_0(L)$ is equal to the number of primes $\ell_i \equiv 1 \pmod p$
dividing $N$ (the genus number of Proposition \ref{r0}). 

\smallskip
Finally, various papers (see \cite{Calegari-Emerton,Lecouturier,
Schaefer-Stubley}) give (for $N \equiv 1 \pmod p$ prime)
a lower bound for $r_1(L)$ as explained in the beginning 
of the Introduction.

\subsubsection{Hubbard--Washington--Schoof results} \label{HWS}

In \cite{HW,Schoof}, some other structural properties of the whole 
class group ${\mathcal H}_M$ of $M$ are given. 

\smallskip
Of course, there are two different cases: the case of the $p$-Sylow submodule 
${\mathcal H}_{M,p} := \Cl_M$ of ${\mathcal H}_M$ and that of the $q$-Sylow 
submodules ${\mathcal H}_{M,q}$, for $q \ne p$. 

\smallskip
Since ${\mathcal H}_K$ 
yields some difficulties, we consider the structure of ${\mathcal H}_M^*$,
the kernel of the algebraic norm $\nu_{M/K} := \sum_{i=0}^{p-1} \sigma^i$; thus
${\mathcal H}_M^*={\mathcal H}_M$ when $K$ is principal ($p=3, 5,...$).
Any intermediate reasonning does exist when some submodules
${\mathcal H}_{K,q}$ of ${\mathcal H}_K$ are trivial. We shall
consider ${\mathcal H}_{M,p}$ in the regular case (for which
${\mathcal H}_{M,p}^*={\mathcal H}_{M,p}=\Cl_M$) and
${\mathcal H}_{M,q}^*$, for all $q \ne p$.

\smallskip
\paragraph{\bf (a) Case $q=p$.}
 This is the non-trivial part which is the purpose of the work 
of Schoof in which the $g$ and $G$-structures are connected by means 
of a ``twisted group ring'' $\Z_p[\zeta_p][g]'$ acting
on the classical filtration of the $p$-class group in $M/K$
(see \cite[\S\,3]{Schoof}). Schoof obtains as corollaries generalizing \cite{HW}:

\smallskip
\quad (i) {\it For $p=3$, if all the $\ell_i$ are inert in $K$, then 
${\mathcal H}_{M,3} \simeq H \times H$
for some abelian $3$-group $H$.}

\smallskip
\quad (ii) {\it For $p=5$, if all the $\ell_i$ are inert in $K$, then there exists
an abelian $5$-group $H$ such that ${\mathcal H}_{M,5} \simeq 
H \times H \times H \times H$, or there exists an exact sequence of the form
$0 \to \F_5 \too {\mathcal H}_{M,5} \too H \times H \times H \times H \to 0$.}

\smallskip
As we shall see, when the primes $\ell_i$ are all inert in $K$, the
filtration of ${\mathcal H}_{M,p} = \Cl_M$ has some more canonical 
behaviour contrary to the general case where the elements of the 
filtration become random. 

\smallskip
For instance, the tables give the following 
more complex examples when some $\ell_i$ split in part:

\footnotesize
\begin{verbatim}
p=3
N=5844=[2,2;3,1;487,1]                 N=7859=[29,1;271,1]
Structure of C_L=[81,3]                Structure of C_L=[81,3]
Structure of C_M=[243,81,3]            Structure of C_M=[81,27,3,3]

p=5 
N=421=[421,1]                          N=451=[11,1;41,1]
Structure of C_L=[5]                   Structure of C_L=[5,5]
Structure of C_M=[5,5,5]               Structure of C_M=[5,5,5,5,5,5]
\end{verbatim}
\normalsize

\begin{remarks}
(i) Assume that $p$ is regular and consider the filtration
$\big(\Cl_i \big)_i$, of $\Cl := \Cl_M$, defined by
$\Cl_{i+1}/\Cl_i :=(\Cl/\Cl_i)^G,\ \, i \geq 0 \ \ \&\ \ \Cl_0=1$.
Then for all $k \geq 1$, the $p^k$-rank $r_k$ (the $\F_p$-dimension of
$\Cl^{p^{k-1}}/\Cl^{p^k}$) is given by
$p^{r_k} = \prd_{i=(k-1)(p-1)}^{k(p-1)-1} \order \big(\Cl_{i+1}/\Cl_i \big)$.

If we express the structure of $\Cl$ under the form:
$$\Cl \simeq \prd_{j=1}^m A/\pi^{n_j}, \ \, 
1\leq n_1 \leq \cdots \leq n_{m}, $$

where $A:=\Z[\zeta_p]$ and $\pi = 1- \zeta_p$, the structure of
abelian group looks like:
$$\Cl \simeq \prd_{j=1}^m \Big[ \big( \Z/p^{a_j+1}\big)^{b_j} 
\oplus \big( \Z/p^{a_j}\big)^{p-1-b_j} \Big], $$

where we have put $n_j = a_j\,(p-1) + b_j$, \ $a_j \geq 0$, 
\ $0 \leq b_j \leq p-2$. 

\smallskip
The $\Cl_i$ appear canonically and some
observations can be done. We have:
$$\Cl_i \simeq \prd_{j,\,n_j \leq i} A/\pi^{n_j} \simeq
 \prd_{j,\,n_j \leq i}\Big[
\big( \Z/p^{a_j+1}\big)^{b_j} \oplus \big( \Z/p^{a_j}\big)^{p-1-b_j}\Big]. $$

Of course, the non-commutative action of $g$ may be precised in 
some cases via the use of the twisted group ring $\Z_p[\zeta_p][g]'$ 
introduced in \cite{Schoof}. 

\smallskip
Using this technique, one may perhaps hope some constraints 
on the elements of the filtration (i.e., on the $a_j, b_j$).

\smallskip
\quad (ii) In the particular case $N=\ell$ prime, $\ell \equiv 1 \pmod p$ 
studied especially in \cite{Calegari-Emerton,Lecouturier,Schoof}, 
we know that ${\rm rk}_p(\Cl_L)\geq 1$ because of $r_0(L)=1$ in 
the relation \eqref{sigma}, and the condition 
${\rm rk}_p(\Cl_L)\geq 2$ depends on $r_1(L)$, where 
$r_1(L) :={\rm dim}_{\F_p}(({\mathcal I}_{M,1})^{e_{\chi_0}})$
with ${\mathcal I}_{M,1} := \Cl_M^{\Theta} \Cl_M^p /
\Cl_M^{\Theta^{2}} \Cl_M^p$. 

\smallskip
More explicit computations and 
conjectures are the purpose of \cite{Lecouturier}.
\end{remarks}

\smallskip
\paragraph{\bf (b) Case $q \ne p$.}
Since the norm $\nu_{M/K}$ commutes with the elements of $g$,
${\mathcal H}_{M,q}^*$ is a $g$-module and a $\Z[G]/(\nu_{M/K})$-module;
so, we can use the structure of $\Z[G]/(\nu_{M/K}) 
\simeq \Z[\zeta_p]$-module as follows. 

\smallskip
Put $A:=\Z[\zeta_p]$ and let $d_q$ be the decomposition group of 
$q$ in $K/\Q$; let ${\mathfrak q}$ be a prime ideal above $q$ in $K$.
We may write:
$$\hbox{${\mathcal H}_{M,q}^* \simeq  \prod_{j =1}^{m}
A/{\mathfrak a}_j$, \ \ ${\mathfrak a}_1 \mid \cdots \mid {\mathfrak a}_m$,} $$

with $q$-ideals ${\mathfrak a}_j$ such that:
$${\mathcal H}_{M,q}^* \simeq  \prd_{j =1}^{m} 
\prd_{s \in g/d_q} A/(s \, {\mathfrak q})^{n_j(s)}, \ \, 
0\leq n_1(s) \leq \cdots \leq n_{m}(s) \hbox{\ for all $s$} . $$
Let ${\mathcal H}$ be a $A$-monogenic factor of 
${\mathcal H}_{M,q}^*$ isomorphic to 
$\prod_{s \in g/d_q} A/(s \, {\mathfrak q})^{n(s)}$.
If ${\mathcal H} = \langle h_0  \rangle_{\Z[G]}$ and
$h = \big(\sum_i a_i \sigma^i \big) \cdot h_0$, then
$t\cdot h = \big (\sum_i a_i\s \sigma^{i \s \omega(t)} \big) 
\cdot ( t \! \cdot\! h_0)$
for all $t \in g$. 

\smallskip
One shows that this twisted operation implies that
the $n_j(s)$ are independent of $s$ (we omit the details; see
\cite[Proposition 2.3]{Schoof}).

\smallskip
Thus:
$$\prd_{s \in g/d_q} A/(s \, {\mathfrak q})^{n_j}
\simeq A/q^{n_j} A \simeq (\Z / q^{n_j}\Z)^{p-1} ; $$ 

whence the result as shown by some examples:

\smallskip
\footnotesize
\begin{verbatim}
N=195=[3, 1; 5, 1; 13, 1] 5 ramified
q=2 ClassqL=[8]
q=2 ClassqM=[8, 8, 8, 8]
q=5 ClassqM=[5]

N=123=[3, 1; 41, 1] 5 ramified
q=3 ClassqL=[9]
q=5 ClassqL=[5, 5]
q=3 ClassqM=[9, 9, 9, 9]
q=5 ClassqM=[5, 5, 5, 5, 5, 5]

N=395=[5, 1; 79, 1] 5 ramified
q=5 ClassqL=[5]
q=7 ClassqL=[7]
q=5 ClassqM=[5, 5]
q=7 ClassqM=[7, 7, 7, 7]

N=623=[7, 1; 89, 1] 5 ramified
q=2 ClassqL=[2]
q=5 ClassqL=[5]
q=11 ClassqL=[11]
q=2 ClassqM=[2, 2, 2, 2]
q=5 ClassqM=[5, 5, 5]
q=11 ClassqM=[11, 11, 11, 11]
\end{verbatim}
\normalsize

\subsection{Reflection theorem}
Since $M$ contains the group of $p$th roots of unity, the
reflection theorem in $M/L$, between the groups $\Cl_M$ and ${\mathcal T}_M$,
is straightforward (\cite[Theorem II.5.4.5]{Gras1}, \cite{Gras4}). 

\smallskip
Recall the correspondence of notations from our book: take $S=\es$, 
$T=P_M$, then $T^*=T$, $S^*= \es$ since the infinite places do not
intervenne for $p>2$, so that this yields:
$$\hbox{$\Cl_{M,T}^S={\rm Gal}(H_M^{\rm pr}/M)$, \ \ 
$\Cl_{M,S^*}^{T^*}=\Cl_M/\cl_M(P_M)$}, $$

and the relation of reflection:
\begin{equation}\label{equa1}
{\rm rk}_{\chi^*}({\rm Gal}(H_M^{\rm pr}/M)) - {\rm rk}_{\chi}(\Cl_M/\cl_M(P_M))
=\rho_{\chi}(P_M, \es),
\end{equation}

where (see \cite[\S\,II.5.4.9.1]{Gras1}) where the base field $k$ is $L$:
\begin{equation}\label{equa2}
\rho_{\chi}(P_M, \es)=\hbox{$\frac{p-1}{2}$}+\sm_{u \in P_{L,\infty}^r} \rho_{u,\chi}
+\sm_{v \in P_{L}} \rho_{v,\chi} +\delta_{\omega,\chi} -  \delta_{\chi_0^{},\chi};
\end{equation}

$\delta$ denotes the Kronecker symbol, $\rho_{w,\chi}=1$
(resp. $0$) if $d_w \subseteq {\rm Ker}(\chi)$ (resp. 
$d_w \not\subseteq {\rm Ker}(\chi)$), $d_w$ being the 
decomposition group of any $w \mid v$ in $M/L$.

\subsubsection{Computation of $\rho_{u,\chi}$ and $\rho_{v,\chi}$}
 In our context, $d_u=\langle c \rangle$ for the 
unique real place $u$ of $L$, where $c$ is the complex conjugation; 
thus we get $\rho_{u,\chi}=1$ (resp. 0) if $\chi$ is even (resp. odd);
in other words:
$$\sm_{u \in P_{L,\infty}^r} \rho_{u,\chi}=\hbox{$\frac{\chi(c)+1}{2}$}. $$

From Kummer congruences on the decomposition 
of $p$ in $M/K$, recall that $p$ is unramified if and only if 
$p \mid N$ or $N^{p-1} \not\equiv 1 \pmod {p^2}$; in this case,
for the unique $p$-place $v$, $d_v=g$, $\rho_{v,\chi}=0$ for any 
$\chi \ne \chi_0^{}$, $\rho_{v,\chi_0}=1$; whence:
$$\sm_{v \in P_{L}} \rho_{v,\chi}=\delta_{\chi_0^{},\chi}. $$

When $p \nmid N$ and $N^{p-1} \equiv 1 \pmod {p^2}$, $N$ is a local $p$th 
power so that $p$ totally split in $M/K$ ($p$ is never inert).
Thus in this case, the decomposition of $p$ in $L$ is of the form:
$$(p)_L = {\mathfrak P}_{L,0} \cdot {\mathfrak P}_{L,1}^{p-1}, $$ 

where ${\mathfrak P}_{L,0}$ is unramified of residue degree $1$ 
and ${\mathfrak P}_{L,1}$ is ramified;
then $d_{v_0}=g$ and $d_{v_1}=1$ giving 
$\rho_{v_0,\chi}=\delta_{\chi_0^{},\chi}$ and $\rho_{v_1,\chi}=1$. 

\medskip
Thus:
$$\sm_{v \in P_{L}}\rho_{v,\chi} = 
\rho_{v_0,\chi} + \rho_{v_1,\chi}=\delta_{\chi_0^{},\chi} + 1. $$

We note that in this case the $p$ completions of $M$ are 
$K_{\mathfrak p}=\Q_p(\mu_p)$, that $L_{v_0}=\Q_p$ and
$L_{v_1}=\Q_p(\mu_p)$.

\subsubsection{Computation of  $\rho_{\chi}(P_M, \es)$}
From \eqref{equa1} and \eqref{equa2} and the above, we get: 

\medskip
(a) $\chi$ even, $p$ totally ramified. Then:
$${\rm rk}_{\chi^*}({\rm Gal}(H_M^{\rm pr}/M)) - {\rm rk}_{\chi}(\Cl_M/\cl_M(P_M))
= \hbox{$\frac{p-1}{2}$} + 1+ \delta_{\omega,\chi}$$ 
 
(a$'$) $\chi$ even, $p$ not totally ramified. Then:
$${\rm rk}_{\chi^*}({\rm Gal}(H_M^{\rm pr}/M)) - {\rm rk}_{\chi}(\Cl_M/\cl_M(P_M))
= \hbox{$\frac{p-1}{2}$} + 2+ \delta_{\omega,\chi}$$  
 
(b) $\chi$ odd, $p$ totally ramified. Then:
$${\rm rk}_{\chi^*}({\rm Gal}(H_M^{\rm pr}/M)) - {\rm rk}_{\chi}(\Cl_M/\cl_M(P_M))
= \hbox{$\frac{p-1}{2}$} + \delta_{\omega,\chi}$$ 

(b$'$) $\chi$ odd, $p$ not totally ramified. Then:
$${\rm rk}_{\chi^*}({\rm Gal}(H_M^{\rm pr}/M)) - {\rm rk}_{\chi}(\Cl_M/\cl_M(P_M))
= \hbox{$\frac{p-1}{2}$} +1+\delta_{\omega,\chi}. $$

\medskip
Which is summarized by the formula:
$${\rm rk}_{\chi^*}({\rm Gal}(H_M^{\rm pr}/M)) - {\rm rk}_{\chi}(\Cl_M/\cl_M(P_M))
= \hbox{$\frac{p-1}{2}$} + d_p+ \hbox{$\frac{\chi(c)+1}{2}$} +\delta_{\omega,\chi}.$$

We note that ${\rm Gal}(H_M^{\rm pr}/M)=\Gamma \oplus {\mathcal T}_M$ as
$g$-modules; thus:
$${\rm rk}_{\chi}({\rm Gal}(H_M^{\rm pr}/M))={\rm rk}_{\chi}({\mathcal T}_M)
+  {\rm rk}_{\chi}(\Gamma_M).$$

We verify that:
$$\hbox{${\rm rk}_{\chi}(\Gamma_M)=\frac{p+1}{2}$\  if $\chi$ is odd,} $$

then:
$$\hbox{${\rm rk}_{\chi}(\Gamma_M)=\frac{p-1}{2}$\  if $\chi \ne \chi_0^{}$ is even,} $$

and:
$$\hbox{${\rm rk}_{\chi_0^{}}(\Gamma_M)=
\hbox{$\frac{p-1}{2}$}+1 = \hbox{$\frac{p+1}{2}$}$.} $$

Thus, for all $\chi$: 
$${\rm rk}_{\chi}(\Gamma_M)=\hbox{$\frac{p-1}{2}$} - \hbox{$\frac{\chi(c)-1}{2}$}
+\delta_{\chi_0^{},\chi}.$$

\begin{theorem}\label{thmrefl}
We have the general relation (for all character $\chi$ of $g$):
$${\rm rk}_{\chi^*}({\mathcal T}_M)=
{\rm rk}_\chi(\Cl_M/\cl_M(P_M)) + d_p  ,$$

where $d_p=1$ (resp. $0$) if $p\nmid N$ and $N^{p-1} \equiv 1 \pmod {p^2}$
(resp. if not).
\end{theorem}

\begin{remark}
If $p \mid N$ or $p \nmid N \ Ê\& \ N^{p-1}\not \equiv 1 \pmod {p^2}$ 
(i.e., $p$ totally ramified in $M/K$), the action of the idempotents
$e_\chi$ is trivial on $\cl_M(P_M)$ for $\chi \ne \chi_0^{}$, giving
${\rm rk}_{\chi^*}({\mathcal T}_M)={\rm rk}_\chi(\Cl_M)$.

\smallskip
For $\chi=\chi_0^{}$, it gives ${\rm rk}_p(\Cl_L/\cl_L(P_L))$;
the class of the prime ideal above $p$ in $L$ is, a priori, 
of order 1 or $p$.

\smallskip
If $p \nmid N$ and $N^{p-1} \equiv 1 \pmod {p^2}$ (i.e., $p$ totally split 
in $M/K$), we have $({\mathfrak p})= {\mathfrak P}_1 \cdots {\mathfrak P}_p$
in $M$ and the $\cl_M({\mathfrak P}_i)$ may be of arbitrary order; in the 
same way, $(p)_L = {\mathfrak P}_{L,0} \cdot {\mathfrak P}_{L,1}^{p-1}$, and
$\cl_L({\mathfrak P}_{L,0})$, $\cl_L({\mathfrak P}_{L,1})$
are unknown.
\end{remark}

\begin{corollary}\label{ranks}
We obtain the following particular cases:
\begin{equation*}
\begin{aligned} 
\hspace{-2.0cm} (i) \hspace{2.3cm} &{\rm rk}_{\omega}({\mathcal T}_M)  =
{\rm rk}_p(\Cl_L/\cl_L(P_L)) + d_p,  \\ 
\hspace{-2.0cm} (ii)\hspace{2.14cm} & {\rm rk}_p({\mathcal T}_L)\  = \,
{\rm rk}_\omega(\Cl_M/\cl_M(P_M))  + d_p, \\
\hspace{-2.0cm}  (iii) \hspace{2cm} &{\rm rk}_p({\mathcal T}_M)  = \,
{\rm rk}_p(\Cl_M/\cl_M(P_M)) +(p-1)\cdot d_p. \\
\end{aligned}
\end{equation*}
\end{corollary}

The relations (ii), (iii), appear in \cite[Proposition III.4.2.2]{Gras1};
indeed, in the case of ${\rm rk}_p({\mathcal T}_L)$, one has
$d_p=\delta_{v_1}=1$ when $N^{p-1} \equiv 1 \pmod {p^2}$
since the ramified completion of $L$ contains $\mu_p$.

\smallskip
We note that for $d_p=1$, ${\mathcal T}_L$ is non-trivial,
which comes from ${\mathcal W}_L \simeq \Z/p\Z$ 
and more generaly we have ${\rm rk}_\chi({\mathcal W}_M)=1$ 
for all $\chi$ (Remark \ref{w} (ii)). 

\smallskip
For instance, for $p=5$ and $d_5=1$, the table gives the following cases:

\footnotesize
\smallskip
\begin{verbatim}
N=124=[2,2;31,1] 5 unramified           N=101=[101,1] 5 unramified
Structure of C_L=[5]                    Structure of C_L=[5]
Structure of T_L=[125]                  Structure of T_L=[25]
Structure of C_M=[5,5]                  Structure of C_M=[5,5,5,5]
Structure of T_M=[125,25,25,25,5,5]     Structure of T_M=[25,25,5,5,5]
\end{verbatim}
\normalsize

\begin{remark}
A number field $F$ is said to be $p$-rational if the conjecture
of Leopoldt holds in $F$ for the prime $p$ and if ${\mathcal T}_F=1$.

\smallskip
See \cite{Gras1, Gras10,Gras3,Gras5,MR} for a full presentation, the 
difficult question of the existence of infinitely many $p$ such that
$F$ is $p$-rational, and an history on this notion which started with 
the pioneer articles \cite{Gras6,Gras7,GJ,JN,Movahhedi,MN}, 
among others.
\end{remark}

\begin{corollary} \label{pprimitive}
If $M$ is $p$-rational then
$L$ and $K$ are $p$-rational\,\footnote{\,Indeed, the transfer map
${\mathcal T}_L  \to {\mathcal T}_M$ is injective for trivial reason 
since $[M:L]$ is prime to $p$, and ${\mathcal T}_K  \to {\mathcal T}_M$ 
is injective under Leopoldt's conjecture in $M$; from \cite{GJ}, 
${\mathcal T}_K=1$ is equivalent to the regularity of $p$.},
$d_p=0$ (i.e., $p$ ramified in $M/K$), then
${\rm rk}_\omega(\Cl_M)=0$, and $\Cl_L=\cl_L(P_L)$.
\end{corollary}

\subsection{Fix points formula for ${\mathcal T}_M$}\label{log}

\subsubsection{Generalities}
We have the following fixed point formula for ${\mathcal T}_F$ in any 
Galois $p$-extension $F/k$ of number fields of Galois group $G$ (under 
the Leopoldt conjecture for $p$ in $F$) \cite[\S\,IV.3, Theorem 3.3]{Gras1}:
\begin{equation}\label{fix}
\order{\mathcal T}_F^G  = \order {\mathcal T}_k \times
\frac{\prod_{{\mathfrak l}\,\notdiv\, p} e^{}_{\mathfrak l}}
{\Big( \sum_{{\mathfrak l}\,\notdiv\, p}\,
\hbox{$\frac{1}{e^{}_{\mathfrak l}}$}  \Z_p {\rm Log}_p({\mathfrak l}) + 
\Z_p {\rm Log}_p(I_k) : \Z_p {\rm Log}_p(I_k) \Big)}, 
\end{equation}

where $e^{}_{\mathfrak l}$ is the ramification index of ${\mathfrak l}$ 
in $F/k$, $I_k$ the group of ideals prime to $p$ of $k$, the 
${\rm Log}_p$-function being defined as follows, where the $k_{\mathfrak p}$
are the $p$-completions of $k$, $E_k$ its unit group and ${\rm log}_p$
the $p$-adic logarithm (\cite[Theorem IV.3.3, Definition IV.3.4]{Gras1}):
$${\rm Log}_p : I_k \too \big(\hbox{$\bigoplus_{{\mathfrak p} \mid p}$} 
k_{\mathfrak p}\big) \big /\Q_p {\rm log}_p(E_k), $$ 

where, for any ideal ${\mathfrak a} \in I_k$ and $m$ such that
${\mathfrak a}^m =: (\alpha)$, $\alpha \in k^\times$, by definition 
${\rm Log}_p({\mathfrak a}) := \frac{1}{m}  {\rm log}_p(\alpha) 
\pmod{\Q_p {\rm log}_p(E_k)}$.
The group ${\mathcal T}_F$ is trivial if and only if ${\mathcal T}_k=1$ and
the condition of ``$p$-primitive ramification'' holds, whence if and only if:
$$\Big(\sm_{{\mathfrak l} \,\nmid\, p} \hbox{$\frac{1}{e_{\mathfrak l}}$} 
{\rm log}_p({\mathfrak l}) + \Z_p {\rm Log}_p(I_k) : \Z_p {\rm Log}_p(I_k) \Big)=
\prd_{{\mathfrak l} \,\nmid\, p}e_{\mathfrak l}. $$

When the base field $k$ is $p$-principal, the condition becomes
${\mathcal T}_k=1$ and:
$$\Big(\sm_{{\mathfrak l} \,\nmid\, p} \hbox{$\frac{1}{e_{\mathfrak l}}$} 
{\rm log}_p({\mathfrak l}) + {\rm Log}_p(U_k) : 
{\rm Log}_p(U_k) \Big)=\prd_{{\mathfrak l} \,\nmid\, p}e_{\mathfrak l}, $$

where $U_k$ is the $\Z_p$-module of principal local units of $k$.

\subsubsection{The case of the extensions $M/K$} \label{M/K}
Note that, for the field $K := \Q(\mu_p)$, ${\rm log}_p(U_K) 
= {\rm log}_p(1 + {\mathfrak p})={\mathfrak p}^2$. We have 
$e_{\mathfrak l}=p$ for each prime ideal ${\mathfrak l} \mid \ell \mid N$.
Thus the condition of $p$-primitive ramification becomes:
$$\Big(\sm_{{\mathfrak l} \,\nmid\, p}\hbox{$\frac{1}{p}$}{\rm log}_p({\mathfrak l}) 
+ {\mathfrak p}^2 + \Q_p{\rm log}_p(E_K) : {\mathfrak p}^2 + 
 \Q_p{\rm log}_p(E_K) \Big)=\prd_{{\mathfrak l} \,\nmid\, p}e_{\mathfrak l}$$

The $\Q_p$-dimension of $\Q_p {\rm log}_p(E_K)$ 
is $\frac{p-1}{2}-1$ (Leopoldt's conjecture holds in $K$) and
${\mathfrak p}^2$ is a free $\Z_p$-module of rank $p-1$.
Thus the $p$-rank of ${\mathcal T}_M$ is non-trivial as soon as
at least $\frac{p-1}{2}+2$ prime ideals ramify in $M/K$.

\subsubsection{Local norms and Hilbert's symbols}\label{symbols}
We shall need the group ${\mathcal N}_{M/K}$ of local norms in $M/K$.
Since $M/K$ is cyclic, ${\mathcal N}_{M/K} = 
{\rm N}_{M/K} (M^\times)$, but the ``product formula''
of norm residue symbols allows us to test local norms
at every place, except the unique $p$-place of $K$ when $p$
ramifies in $M/K$.

\smallskip
Recall that the norm residue symbol in $M/K$, at a place $v$, is the map:
$$(\ \bullet\  , K_v^\times (\sqrt[p]{N})/K_v) : K_v^\times \too {\rm Gal}(M_v/K_v)$$ 

for which the image of $y \in K^\times$ is trivial if and only if 
$y$ is local norm at $v$.
Since $K$ contains $\mu_p$, this symbol is characterized by the
Hilbert's symbol
$(\ \bullet\  , \ \bullet\  ) : K_v^\times \times K_v^\times \too \mu_p$,
defined by:
$$(x, y)_v := \frac{(y ,\  K_v(\sqrt[p]{x})/K_v) \cdot \sqrt[p]{x} }{\sqrt[p]{x}} \in \mu_p.$$

The main properties of these symbols are well known 
(see, e.g., \cite[\S\,II.7.1]{Gras1}); for instance, for $x=N$,
consider $(N, y)_v$ 
when $v$ is the prime ideal ${\mathfrak l}$ of $K$
of residue field $\F_{\ell^f}$, unramified in 
$K_{\mathfrak l}(\sqrt[p]{y})/K_{\mathfrak l}$ 
and with ${\rm v}_{\mathfrak l}(N) \not\equiv 0 \pmod p$ (e.g., 
${\mathfrak l} \mid \ell \mid N$); then the symbol 
depends on the Frobenius in the local extension and is trivial 
if and only if this Frobenius is trivial, thus $y$ local $p$th power in 
$K_{\mathfrak l}^\times$; in other words if and only if
$y^{\frac{\ell^f-1}{p}} \equiv 1 \pmod {\mathfrak l}$.

\smallskip
We will also use the ``conjugation property'':
$$(x,\, y)_{\mathfrak l}^s = (x^s,\, y^s)_{s\, \cdot\, {\mathfrak l}}, $$ 

for any isomorphism $s$ of $K$.

\subsubsection{Criteria of $p$-rationality and $p$-principality of $M$}
We assume that the Leopoldt conjecture holds for $p$ in $M=K(\sqrt[p]{N})$.

\begin{theorem} \label{equiv}
Let $E_K^{{\mathfrak p}} = E_K \cdot \langle 1-\zeta_p  \rangle_\Z$ 
be the group of $p$-units of $K$.

\smallskip
({\bf a}) The following conditions are equivalent:

\smallskip
\quad (i) ${\mathcal T}_M =1$ (i.e., $M$ is $p$-rational).

\smallskip
\quad (ii) ${\mathcal T}_K =1$ and $M/K$ is $p$-primitively ramified.

\smallskip
\quad (iii) $p$ is regular and $M/K$ is $p$-primitively ramified.

\smallskip
\quad (iv) $p$ is ramified in $M/K$ and $\Cl_M = \cl_M({\mathfrak P}_M)$
where ${\mathfrak P}_M \mid p$ in $M$.

\smallskip
\quad (v) $M/K$ and $p$ fulfill the following conditions:

\smallskip
\qquad \ $\bullet$ $p$ is regular,

\smallskip
\qquad \ $\bullet$ $p$ is ramified in $M/K$, 

\smallskip
\qquad \ $\bullet$ $(E_K^{{\mathfrak p}} : E_K^{{\mathfrak p}} \cap {\mathcal N}_{M/K})
=\prod_{{\mathfrak l} \, \nmid \, p} e_{\mathfrak l}$.

\smallskip
({\bf b}) We have ${\mathcal T}_M = \Cl_M =1$ if and only if all the following
conditions hold: 

\quad (i) $p$ is regular,

\smallskip
\quad (ii) $p$ ramifies in $M/K$,

\smallskip
\quad (iii) $(E_K^{{\mathfrak p}} : E_K^{{\mathfrak p}} \cap {\mathcal N}_{M/K})=
(E_K : E_K \cap {\mathcal N}_{M/K}) =
\hbox{$\prod_{{\mathfrak l} \, \nmid \, p}$} e_{\mathfrak l}$ and
this occurs if and only if the two following conditions are fulfilled:\par
\qquad \ $\bullet$ 
There exists $\varepsilon_0 \in E_K$ such that 
$(p\cdot \varepsilon_0)^{\frac{\ell^{f_{\mathfrak l}}-1}{p}} \equiv 1 
\pmod {{\mathfrak l}}$, for all ${\mathfrak l} \nmid p$, ramified in $M/K$,
where $f_{\mathfrak l}$ is the residue degree of ${\mathfrak l}$ in $K/\Q$;

\smallskip
\qquad \ $\bullet$ 
$(E_K : E_K \cap {\mathcal N}_{M/K}) =
\hbox{$\prod_{{\mathfrak l} \, \nmid \, p}$} e_{\mathfrak l}$.
\end{theorem}

\begin{proof}  (a) We have $E_K^{{\mathfrak p}}= E_K \cdot \langle p \rangle$,
up to $p$-powers, since $(1-\zeta_p)^{p-1} = p \cdot \eta$, $\eta \in E_K$. 

\smallskip
We have (i) $\Longleftrightarrow$ (ii) (Corollary \ref{pprimitive}); then
(ii) $\Longleftrightarrow$ (iii) since in this case the $p$-rationality is
equivalent to the $p$-regularity \cite{GJ}; from
Corollary \ref{ranks} to Theorem \ref{thmrefl}, we get
(i) $\Longleftrightarrow$ (iv);
for (i) $\Longleftrightarrow$ (v) see \cite[Proposition IV.4.8.2]{Gras1}.

(b)  Assume ${\mathcal T}_M = \Cl_M =1$.
The regularity of $p$ is clear, what we assume.
From the Corollary \ref{ranks}, ${\mathcal T}_M=1$
implies $d_p=0$ ($p$ totally ramified); thus (Chevalley's formula):
$$\ds \order \Cl_M^G = \order \Cl_K \!\cdot\! \frac{p \!\cdot\! 
\prod_{{\mathfrak l} \, \nmid \, p} e_{\mathfrak l}}
{p\!\cdot\! (E_K : E_K \cap {\mathcal N}_{M/K})} = 1; $$

whence $(E_K : E_K \cap {\mathcal N}_{M/K}) =
\hbox{$\prod_{{\mathfrak l} \, \nmid \, p}$} e_{\mathfrak l}$.
The map:
$$E_K/E_K \cap {\mathcal N}_{M/K} \too 
E_K^{{\mathfrak p}} / E_K^{{\mathfrak p}} \cap {\mathcal N}_{M/K}$$

is injective; thus the left equality of the statement holds if and only if:
$$E_K^{{\mathfrak p}} = 
E_K \!\cdot\! (E_K^{{\mathfrak p}} \cap {\mathcal N}_{M/K}). $$

Since $E_K^{{\mathfrak p}}=E_K \cdot \langle p \rangle$ (up to $p$-powers), 
the equality holds if and only if there exists $\varepsilon_0 \in E_K$ 
such that $p\!\cdot\! \varepsilon_0 \in {\mathcal N}_{M/K}$.

\smallskip
Since there exists a unique $p$-place, from the product formula
the condition becomes that there exists $\varepsilon_0 \in E_K$ 
such that (in terms of Hilbert's symbols),
$\big(N, p\!\cdot\! \varepsilon_0 \big)_{\mathfrak l}=1$, 
for all ${\mathfrak l} \nmid p$, ramified in $M/K$;
in other words:
$$(p\!\cdot\! \varepsilon_0)^{\frac{\ell^{f_{\mathfrak l}}-1}{p}} 
\equiv 1 \pmod {{\mathfrak l}},$$

for all ${\mathfrak l} \nmid p$ ramified in $M/K$.
\end{proof}

\subsection{The $3$-rationality of $L=\Q(\sqrt[3]{N})$}
We intend to characterize the integers $N$ such that $L$ is
$3$-rational and show that ${\mathcal T}_L=1$ implies 
${\mathcal T}_M=1$. So, we assume ${\mathcal T}_L=1$.

\subsubsection{Generalities}

Let $g = {\rm Gal}(M/L) = \{1, c\}$, where $c$ is the complex conjugation;
thus $e_{\chi_0^{}} = \frac{1+c}{2}$ and $e_{\omega} = \frac{1-c}{2}$. 

\smallskip
If ${\mathcal T}_L=1$, Corollary \eqref{ranks}\,(ii) implies
the two following properties:

\smallskip
\quad (i) $d_3=0$ (i.e., $3$ is ramified in $M/K$), which implies $N=3\,N'$,
or $9\,N'$, $N'$ prime to $3$, or $N=\pm1 + 3\,u$, $3 \nmid u$, 

\smallskip
\quad (ii) $\Cl_M^{e_\omega}=1$ (indeed, $d_3=0$ implies
$P_M = \{{\mathfrak P}_M\}$ 
and $\cl_M({\mathfrak P}_M)^{e_\omega}=1$ since 
${\mathfrak P}_M$ is invariant under $g$).

\smallskip
Thus $\Cl_M = \Cl_M^{e_{\chi_0}} \oplus
\Cl_M^{e_{\omega}} = \Cl_M^{e_{\chi_0}} \simeq  \Cl_L$.
Moreover, the $3$-Hilbert class field $H_L$ of $L$ must be contained in
the compositum $\wt L$ of its two independent $\Z_3$-extensions 
(from the Schema in Lemma \ref{exact}, $\ov \Cl_L \ne 1$ implies 
${\mathcal T}_L\ne1$).

\begin{lemma}\label{lem20}
The cyclotomic $\Z_3$-extension $M_\infty$ of $M$ is 
totally ramified at~$3$; whence $H_M \cap M_\infty = M$, 
$H_L \cap L_\infty = L$ and 
${\rm rk}_3(\Cl_M) = {\rm rk}_3(\Cl_L) \leq 1$.
\end{lemma}

\begin{proof}
Consider the cyclic cubic extension $\Q_1$ of $\Q$ in $\Q_\infty$
and the extensions $L\Q_1/\Q$ and $M\Q_1/K$.
For $3 \nmid N$, the radicals $\zeta_3, N, N \zeta_3, N\zeta_3^2$
in $K$ all give ramified extensions at $3$ (Kummer congruences); 
the case $3 \mid N$ being obvious.
So the extension $M(\sqrt[3]{\zeta_3})=M \Q_1$ is ramified over $M$,
otherwise all the four relative cubic extensions in $M\Q_1/K$
should be unramified, which is absurd because of the existence
of a non-trivial inertia group. Thus this goes down to $L\Q_1/L$.\end{proof}

\subsubsection{Technical lemmas}\label{technical}
Under the assumption ${\mathcal T}_L=1$,
one has obtained ${\rm rk}_3(\Cl_M) = {\rm rk}_3(\Cl_L) \leq 1$, 
and the ramification of $3$ in $M/K$.

\begin{lemma}\label{lem22} 
The $3$th root $\zeta_3$ is norm in $M/K$ if and only if 
all $\ell \mid N$ are such that $\ell \equiv \pm1 \pmod 9$, 
which implies necessarily $3 \mid N$.

\smallskip
The case $\Cl_M=1$ is equivalent to 
$\prod_{{\mathfrak l}\, \nmid \, 3} e_{\mathfrak l} = 
(\langle \zeta_3 \rangle : \langle \zeta_3 \rangle
\cap {\mathcal N}_{M/K})$ (Chevalley's formula); 
so $\zeta_3$ norm gives $N=3$, otherwise, 
$\zeta_3$ non-norm gives 
$N=\ell$, $3\, \ell$ or $9\, \ell$, 
with $\ell = -1+ 3\, u$, $3 \nmid u$.
When $\Cl_M \ne 1$, necessarily
$\Cl_M= \Cl_M^G \simeq \Z/3\Z$ and
$\prod_{{\mathfrak l}\, \nmid \, 3} e_{\mathfrak l}
= 3\cdot (\langle \zeta_3 \rangle : \langle \zeta_3 \rangle
\cap {\mathcal N}_{M/K}) \in \{3,9\}$.
\end{lemma}

\begin{proof}
The first claims are obvious since $3$ must be ramified in $M/K$
Assume now $\Cl_M) \ne 1$, thus ${\rm rk}_3(\Cl_M) = 1$.
Since ${\rm N}_{M/K}(\Cl_M)=1$, we consider $\Cl_M$ as a
$\Z_3[\zeta_3]$-module and write:
$$\Cl_M \simeq \prd_{j =1}^m
\Z_3[\zeta_3]/\pi^{n_j}\Z_3[\zeta_3], \ \, 
1\leq n_1 \leq \cdots \leq n_m , $$ 

where $\pi = 1-\zeta_3$; but to get the $3$-rank
$1$, we must have $m=1$ and $n_1=1$, whence
$\Cl_M=\Cl_M^G \simeq \Z_3[\zeta_3]/\pi\Z_3[\zeta_3] \simeq \Z/3\Z$.
\end{proof}

\begin{lemma} In the cases
$\prod_{{\mathfrak l}\, \nmid \, 3} e_{\mathfrak l} = 9$, 
we have $\Cl_M^G = \cl_M(I_M^G) \simeq \Z/3\Z$, where 
$I_M^G$ is the group of $G$-invariant ideals of $M$.
\end{lemma}

\begin{proof} 
The condition $\prod_{{\mathfrak l}\, \nmid \, 3} e_{\mathfrak l} = 9$
is equivalent to $\zeta_3$ non-norm. 
From the classical exact sequence:
$$1\too \cl_M(I_M^G) \too \Cl_M^G \too 
E_K \cap {\mathcal N}_{M/K} \big / {\rm N}_{M/K}(E_M) \too 1,$$

we deduce the claim since $E_K \cap {\mathcal N}_{M/K} =
\langle \zeta_3 \rangle \cap {\mathcal N}_{M/K} =1$.
\end{proof}

\begin{lemma}\label{lem2} The case $N=3 \cdot N'$, $N' \ne 1$,
with $N' \equiv \pm 1 \pmod 9$, implies $\Cl_M^{e_\omega} \ne 1$
and is then to be excluded.
\end{lemma}

\begin{proof} One verifies that the extension $M(\sqrt[3]{3})/M$
is an unramified extension of degree $3$ (the subextension
$K(\sqrt{N'})/K$ is unramified at $3$ and $K(\sqrt[3]{3})/K$ is
unramified outside $3$); so $M(\sqrt[3]{3})=H_M$. 
Since the character of the radical $\langle 3 \rangle M^{\times 3}$ 
is $\chi_0$, that of ${\rm Gal}(M(\sqrt[3]{3})/M) \simeq \Cl_M$ is 
$\omega$ (absurd).
\end{proof}

The case $\prod_{{\mathfrak l}\, \nmid \, 3} e_{\mathfrak l} = 3$,
when $\zeta_3$ is non-norm, is related to the first case $\Cl_M=1$
for which we have only 
$N= \ell$, $3\,\ell$ or $9\,\ell$, with $\ell = -1 + 3u$, $3 \nmid u$. 
The second case $\prod_{{\mathfrak l}\, \nmid \, 3} e_{\mathfrak l} = 3$ 
and $\zeta_3$ norm is impossible  (take into account the ramification 
of $3$ for $N=\ell$ giving $3 \nmid u$, and Lemma \ref{lem2} 
for $3 \mid N$).

\smallskip
We shall see in the forthcomming \S\,\ref{364} that in the case 
$\prod_{{\mathfrak l}\, \nmid \, 3} e_{\mathfrak l} = 9$ with
$N=3^\delta \ell$, $\ell \equiv 1 \pmod 3$,
some supplementary conditions are necessary.

\smallskip
We examine now the remaining case for
$\prod_{{\mathfrak l}\, \nmid \, 3} e_{\mathfrak l} = 9$
where two prime divisors of $N$ are inert in $K/\Q$, 
thus we have $N = \ell_1^{\lambda_1} \ell_2^{\lambda_2}$,
$3\,\ell_1^{\lambda_1} \ell_2^{\lambda_2}$ or 
$9\,\ell_1^{\lambda_1} \ell_2^{\lambda_2}$, with
$\ell_i \equiv -1 \pmod 3$, $\lambda_i \not\equiv 0 \pmod 3$.

\smallskip
Consider the filtration $\big(\Cl_i =: \cl_M({\mathcal I}_i) \big)_i$, 
of $\Cl := \Cl_M$, where ${\mathcal I}_i$ is a suitable 
sub-module of $I_M$, and defined by:
$$\Cl_{i+1}/\Cl_i :=(\Cl/\Cl_i)^G,\ \, i \geq 0 \ \ \&\ \ \Cl_0=1.$$

Then, under the ramification of $3$, we have the formula:
$$\order (\Cl_{i+1}/\Cl_i) = \ds
\frac{\prod_{{\mathfrak l}\, \nmid \, 3} e_{\mathfrak l}}
{(\Lambda_i : \Lambda_i \cap {\mathcal N}_{M/K})}, \ \ 
\Lambda_i := \{x \in K^\times,\, (x) \in {\rm N}_{M/K}({\mathcal I}_i)\}, $$
with ${\mathcal I}_0=1$, $\Lambda_0=E_K$; then ${\mathcal I}_1$ 
is the set of representatives of $\Cl_1$ (it contains $I_M^G$)
allowing the determination 
of $\Lambda_1$, etc. (see \cite{Gras8}). The equalities $\Cl_M=
\Cl_M^G \simeq \Z/3\Z$ means that the ``algorithm'' must stop at $\Cl_2$
(i.e., $\Cl_2=\Cl_1=\Cl$), a condition that we will use below. 

\smallskip
So $\Cl_M = \langle \cl_M(I_M^G) \rangle$ and
$I_M^G$ is generated by the prime ideal above $3$, 
the $\ell_i \mid N$ and the extension
of the ideals of $K$ (which are $3$-principal but do not intervenne):
$$I_M^G = \langle {\mathfrak P}, {\mathfrak L}_1, {\mathfrak L}_2 \rangle 
\  (I_K). $$

Then ${\rm N}_{M/K}(I_M^G) = \langle {\mathfrak p}, (\ell_1), (\ell_2) \rangle$,
up to $I_K^3$, giving (since ${\mathfrak p}^2=(3)$):
$$\Lambda_1 = \langle \zeta_3, 3,  \ell_1, \ell_2 \rangle \cdot K^{\times 3}.$$

We consider the Hilbert symbols $(N,3)_{\ell_j}$ and
$(N,\ell_i)_{\ell_j}$, $i, j \in \{1,2\}$. 

\smallskip
If for instance $(N,\ell_i)_{\ell_j} = \zeta'_3$,
by conjugation, one gets $(N,\ell_i)_{\ell_j} = \zeta'_3{}^c = \zeta'_3{}^{-1}$ 
since $\ell_i, N, \ell_j$ are invariants by $c$ (see \S\,\ref{symbols}); 
so all these symbols are trivial. 

So $(\Lambda_1 : \Lambda_1 \cap {\mathcal N}_{M/K})=3$,
because of $\zeta_3$, and 
$\ds\frac{\prod_{\mathfrak l \, \nmid \, 3} e_{\mathfrak l}}
{(\Lambda_1 : \Lambda_1  \cap {\mathcal N}_{M/K})} = \frac{9}{3} =3$,
which means that $\Cl_2/\Cl_1$ is of order $3$ and that 
$\Cl_M$ is non-cyclic (absurd).

\smallskip
So, we may now, for proving the property
${\mathcal T}_L=1$ implies ${\mathcal T}_M=1$, restrict 
ourselves to the cases $N= \ell$, $3\,\ell$ or $9\,\ell$, with 
$\ell \equiv \pm 1 +3\,u$, $3 \nmid u$.

\subsubsection{Case 
$N=3^{\delta} \ell$, with $\ell = -1 + 3u$, $3 \nmid u$}
At this step, recall that this case implies 
${\mathcal T}_M=1$; since $\zeta_3 \notin  {\mathcal N}_{M/K}$,
we get $\Cl_M^G = 1$, thus $\Cl_M = 1$ giving the result
(Corollary \ref{ranks}\,(iii)).
We shall give the direct characterization of the $3$-rationality
of $M$, using the fixed point formula \eqref{fix} for
${\mathcal T}_M$ in $M/K$, and note that the condition gives
exactly the above values of $N$ when $\ell$ is inert in $K$; 
so this formula will also give the triviality of ${\mathcal T}_M$.

\subsubsection{Case $N=3^{\delta} \ell$, with 
$\ell = 1 + 3 \, u$, $3 \nmid u$}\label{364}

In this case, we will see that one must eliminate some 
values of $\ell$; indeed,
for this family, $d_3=0$ (ramification of $3$ in $M/K$),
$\zeta_3 \notin {\mathcal N}_{M/K}$
and the table in the Appendix shows that the $3$-rationality 
holds in many cases (e.g., for $N= 3, 7, 13, 19, 31,\,\ldots$\,), 
but not for $N=61=1+ 3 \cdot 20 = 
{\rm N}_{K/\Q}\big (\frac{1+ 9\,\sqrt{-3}}{2} \big)$. 

\begin{lemma}\label{lem4} Put $\ell = \alpha \cdot \alpha'$, with
$\alpha =A+B\sqrt{-3}$, where $A, B$ are integers or half integers.
Then considering $\zeta_6\! \cdot \!\alpha$, $\zeta_6 \in \mu_6$,
we may suppose that $A=1+3\,a$ and $B=3\,b$, $a, b \in \Z$ or 
$\frac{1}{2} \Z$.
\end{lemma}

\begin{proof} Assume that $\pm \alpha$ is not of the above form and
write $\zeta_6 = s \frac{-1+ s' \sqrt{-3}}{2}$ where $s,s' \in \{\pm1\}$; then:

\centerline{$\zeta_6 \cdot \alpha = 
\ds s \Big[ \frac{-A-3\,s' B + (s'A - B) \sqrt{-3}}{2}  \Big].$}

\smallskip
Since $A$ and $B$ are prime to $3$, there exists $s'$ such that
$s'A - B \equiv 0 \pmod 3$. The existence of $s$ such that
$s \, \frac{-A-3\, s' B}{2} \equiv 1\pmod 3$ is clear.
\end{proof}

So we assume that $\alpha = 1+3\, a+ 3\, b \sqrt{-3}$
where $a, b$ are integers or half integers. We then have
$\ell = (1+3\,a)^2+27\, b^2 \equiv 1+6\,a \pmod 9$.
So, if $N=\ell$, the condition of ramification of $3$ becomes
$a \not \equiv 0 \pmod 3$. The case $3 \mid N$ 
is given by the Lemma \ref{lem2}; so, in any case, 
$a \not \equiv 0 \pmod 3$.

\begin{lemma}\label{lem5} The case $b \equiv 0 \pmod 3$
must be excluded.
\end{lemma} 

\begin{proof}
If $b \equiv 0 \pmod 3$, then 
$\alpha \cdot \alpha'{}^{-1} \equiv 1 \pmod {9\, {\mathfrak p}}$ in $K$, giving
$\alpha \cdot \alpha'{}^{-1} \equiv 1 \pmod {{\mathfrak P}^{15}}$ in $M$
with ${\mathfrak P} := {\mathfrak P}_M \mid 3$; then
we see that $(\alpha \cdot \alpha'{}^{-1})$ is the cube of the ideal
${\mathfrak L}\,{\mathfrak L}'^{-1}$ of $M$ , ${\mathfrak L} \mid {\mathfrak l}$,
${\mathfrak L}' \mid {\mathfrak l}'$,
and that $\alpha \cdot \alpha'{}^{-1} \notin M^{\times 3}$; indeed,
$\alpha \cdot \alpha'{}^{-1} = x^3$ implies
$x = \sqrt[3]{N}^\lambda \cdot \beta$, $\beta \in K^\times$, which
gives ${\mathfrak l}\,{\mathfrak l}'^{-1} = (N^\lambda \beta^3) =
{\mathfrak l}^\lambda\,{\mathfrak l}'^{\lambda} \beta^3$ in $K$,
which implies $\lambda \equiv 1 \pmod 3$ and 
$\lambda+1 \equiv 0 \pmod 3$ (absurd).

\smallskip
So, as radical in $M$, $\alpha \cdot \alpha'{}^{-1}$ defines 
the unramified cubic cyclic extension$H_M$ of $M$ in which $3$ 
splits (indeed, for the non-ramification of $3$ one needs the 
congruence modulo ${\mathfrak P}^{9}$ and the splitting is 
then obvious); since $\Cl_M = \Cl_M^G \simeq \Z/3\Z$
from Lemma \ref{lem22}, 
the ideal ${\mathfrak P}$ (split in $H_M$) is principal, and 
$\cl_M(I_M^G) = \langle \cl_M({\mathfrak L}), \cl_M({\mathfrak L}^c ) \rangle)
\simeq \Z/3\Z$; but the relations ${\mathfrak L}\,{\mathfrak L}^c = (\sqrt[3]{\ell})$ 
or ${\mathfrak P}^{\delta}{\mathfrak L}\,{\mathfrak L}^c =
(\sqrt[3]{3^{\delta}\ell})$ show that $\cl_M({\mathfrak L}\,{\mathfrak L}^c )=1$
with $\cl_M({\mathfrak L})$ of order $3$, giving a non trivial class in 
$\Cl_M^{e_\omega}$ (absurd).
\end{proof}

\subsubsection{Use of the fixed points formula}
Since the ``good candidates'' obtained by elimination for our 
purpose will be given by the direct forthcomming computation 
of $\order {\mathcal T}_M^G$, we may state:

\begin{theorem}\label{3rat}
For $p=3$, the fiield $L=\Q(\sqrt[3]{N})$ is $3$-rational if and only 
if $N$ fulfills one of the three following possibilities (up to $\Q^{\times 3}$):

\smallskip
(ii) $N=3$,

\smallskip
(ii) $N=\ell$, $3\, \ell$ or $9\, \ell$, with $\ell = -1+ 3\, u$, $3 \nmid u$,

\smallskip
(iii) $N=\ell$, $3\, \ell$ or $9\, \ell$, with $\ell = 
{\rm N}_{\Q(\sqrt{-3})/\Q}\big (1+3\,( a+b \sqrt{-3}) \big)$, 
with $3 \nmid a \s b$.

\smallskip
Moreover, this occurs if and only if $M=L(\sqrt{-3})$ is $3$-rational. 
\end{theorem}

\begin{proof}
Let ${\mathcal O}_K$ be the ring of integers of the $3$-completion 
$K_{\mathfrak p}$ of $K$ (see \S\,\ref{M/K}); the condition 
of $3$-rationality of $M$ (i.e., of $3$-primitive ramification) becomes:
$$\Big(\sm_{{\mathfrak l} \,\nmid\, 3} \hbox{$\frac{1}{3}$} {\rm log}({\mathfrak l})
+3\, {\mathcal O}_K : 3\, {\mathcal O}_K \Big)=3^t, \ t=t_1+ 2 \,t_2, $$

where $t_1$ (resp. $t_2$) is the number of $\ell \equiv -1 \pmod 3$
(resp. $\ell \equiv 1 \pmod 3$) dividing $N$.
We have $t \leq 2$ since the denominator is at most $9$
because ${\rm log}({\mathfrak l}) \in {\mathcal O}_K$; 
thus $N$ is of one of the following forms (besides $N=3$):

\smallskip
(i) $N \in \{\ell, \, 3\,\ell ,\,9\,\ell , \ \, \ell \equiv -1 \pmod 3 \}$,

\smallskip
(ii) $N \in \{\ell_1^{\lambda_1}\ell_2^{\lambda_2}, \, 
3\, \ell_1^{\lambda_1}\ell_2^{\lambda_2},\, 
9\,\ell_1^{\lambda_1}\ell_2^{\lambda_2}, \ \, \ell_i \equiv -1 \pmod 3,\ 
\lambda_i \not\equiv 0 \pmod 3 \}$,

\smallskip
(iii)  $N \in \{\ell, \, 3\,\ell ,\,9\,\ell , \ \, \ell \equiv 1 \pmod 3 \}$.

\medskip
Case (i). The prime $\ell$ is inert in $K$, thus
$\prod_{{\mathfrak l} \,\nmid\, 3}e_{\mathfrak l}=3$, and the condition
${\mathcal T}_M^G=1$ holds if and only if
${\rm log}(\ell) \in 3\, \Z_3^\times$ to get
$\big ( \frac{1}{3} {\rm log}(\ell)+3\, {\mathcal O}_K : 3\, {\mathcal O}_K \big)=3$,
hence $\ell = -1+ 3\, u$, $3 \nmid u$.

\medskip
Case (ii). The previous computation shows that this case (where 
$\prod_{{\mathfrak l} \,\nmid\, 3}e_{\mathfrak l}=9$) 
gives the non-$3$-rationality
of $M$ since the logarithms take values only in $3\, \Z_3$.

\medskip
Case (iii).  Here $\ell$ splits which yields 
$\prod_{{\mathfrak l} \,\nmid\, 3}e_{\mathfrak l}=9$.
Put $\frac{1}{3}{\rm log}({\mathfrak l}) = u + v \sqrt{-3}$ in
${\mathcal O}_K=\Z_3\oplus \Z_3 \sqrt{-3}$; so 
$\frac{1}{3}{\rm log}({\mathfrak l}^c) = u - v \sqrt{-3}$. 
To get the denominator $9$ we must have:
$$\big\langle u + v \sqrt{-3}, u - v \sqrt{-3} \big\rangle_{\Z_3}^{} =
\Z_3\oplus \Z_3 \sqrt{-3}, $$ 

thus $u$ and $v$ prime to $3$.
Put ${\mathfrak l} = (\alpha_\ell)$ with (Lemma \ref{lem4}):
$$\hbox{$\alpha_\ell = 1+3\,( a+b \sqrt{-3})$, \ $a, b \in \frac{1}{2}\Z$}; $$

giving $u \equiv a \!\pmod 3$, $v \equiv b\! \pmod 3$ and finds
again the condition on $a$, $b$; whence the result and the equivalence
${\mathcal T}_L=1$ if and only if ${\mathcal T}_M=1$.
\end{proof}

The following program computes the previous parameters $a, b$, and
allows the verifcation of the criterion when ${\mathcal T}_L=1$.

\footnotesize
\begin{verbatim}
{PK=x^2+3;K=bnfinit(PK,1);z6=(1+Mod(x,PK))/2;n=12;
forprime(N=2,10^3,if(Mod(N,9)!=4 & Mod(N,9)!=7,next);PL=x^3-N;
L=bnfinit(PL,1);Lpn=bnrinit(L,3^n);H=Lpn.cyc;T=List;e=matsize(H)[2];
for(k=1,e-2,c=H[e-k+1];if(Mod(c,3)==0,listinsert(T,3^valuation(c,3),1)));
nu=bnfisintnorm(K,N)[1];nu=Mod(nu,PK);
for(k=1,6,nu=z6*nu;A=component(lift(nu),1);B=component(lift(nu),2);
if(Mod(A,3)!=1||Mod(B,3)!=0,next);a=Mod((A-1)/3,3);b=Mod(B/3,3);
if(T!=List([]),
print("N=",N,"  a=",lift(a),"  b=",lift(b),"   L non-3-rational"));
if(T==List([]),
print("N=",N,"  a=",lift(a),"  b=",lift(b),"   L 3-rational"));break))}
\end{verbatim}
\normalsize

\subsection{The $3$-principality of $L=\Q(\sqrt[3]{N})$}
After the criterion of $3$-rationa\-lity, we will prove an analogous
result for the triviality of $\Cl_L$, $L = \Q(\sqrt[3]{N})$.

\begin{theorem} \label{class}
Let $N=3^\delta \prod_{i=1}^n \ell_i^{\lambda_i}$, 
$\delta \geq 0$, $n \geq 0$,
$\lambda_i \not\equiv 0 \pmod 3$.
The $3$-class group of $L = \Q(\sqrt[3]{N})$ is trivial 
if and only if that of its Galois closure $M$ is trivial,
which is equivalent to $\ell_i\equiv -1 \pmod 3$ for all $i$ and
$n = \nu_3 + d_3$, where $\nu_3=0$ (resp. $\nu_3=1$) 
if $\ell_i \equiv -1 \pmod 9$ for all $i$ (resp. if not), where 
$d_3=1$ (resp. $d_3=0$) if $N \equiv \pm1 \pmod 9$ (resp. if not).
\end{theorem}

\begin{proof} Assume that $\Cl_L=1$ and consider $\Cl_M$;
from Proposition \ref{r0}, all the $\ell_i$ 
are such that $\ell_i \equiv -1 \pmod 3$.
We then have $\Cl_M=\Cl_M^\omega$ by assumption.

\smallskip
The Jaulent exact sequence \eqref{jfj} yields
$(\Cl_M^G)^\omega = \Cl_M^\omega$ and $(\Cl_M^\omega)^\Theta  =1$; 
whence $\Cl_M=\Cl_M^G$, since $\Cl_M=\Cl_M^\omega$
and $\Cl_M^\omega \subseteq \Cl_M^G$ (i.e., $\Cl_2=\Cl_1$ in 
terms of the filtration of $\Cl_M$).

\smallskip
Chevalley's formula gives
$\order \Cl_M^G = 3^{n-d_3-\nu_3}$ since
$(\langle \zeta_3 \rangle\! :\! \langle \zeta_3 \rangle 
\cap {\mathcal N}_{M/K})=3^{\nu_3}$ and $d_3=0$
(resp. $1$) if $3$ ramifies (resp. if not).

\smallskip
So, assuming $\Cl_M \ne 1$ implies $\Cl_M^G \ne 1$, whence:
$$n-d_3 \geq \nu_3+1. $$

Moreover, $\nu_3=0$ implies $d_3=1$ (Lemma \ref{lem22}).
Consider two cases:

\smallskip
\quad (i) Case $\nu_3=1$ ($\zeta_3 \notin {\mathcal N}_{M/K}$, 
$d_3 \in \{0,1\}$ and $n-d_3 \geq 2$).
The filtration leads to 
$\Lambda_1 = \langle \zeta_3, \ell_1, \ldots , \ell_n  \rangle$ if $d_3=1$,
$\Lambda_1 = \langle \zeta_3, 3, \ell_1, \ldots, \ell_n  \rangle$ if not; so:
$$\order (\Cl_2/\Cl_1) = \ds \frac{3^{n-d_3}}
{(\Lambda_1 : \Lambda_1 \cap {\mathcal N}_{M/K})}. $$

We obtain that
$(\Lambda_1 : \Lambda_1 \cap {\mathcal N}_{M/K})=3$,
whatever the ramification of $3$, since
all Hilbert's symbols with rational integers, at inert places, are trivial;
so $\order (\Cl_2/\Cl_1) = 3^{n-d_3-1}$.
Since $n-d_3 -1\geq \nu_3 =1$, we obtain
$\order (\Cl_2/\Cl_1) \equiv 0 \pmod 3$ (absurd).

\smallskip
\quad (ii) Case $\nu_3=0$ ($\zeta_3 \in {\mathcal N}_{M/K}$, 
$d_3=1$ and $n\geq 2$).
The relation $\zeta_3 = {\rm N}_{M/K}(y)$,
$y \in M^\times$, leads to $(y) = {\mathfrak A}_M^{1-\sigma}$, for
an ideal ${\mathfrak A}_M$ of $M$. Since 
$\zeta_3=\zeta_3^{e_\omega}$, we may write 
$(y^{e_\omega}) = {\mathfrak A}_M^{(e_\omega(1-\sigma))} = 
{\mathfrak A}_M'^{\s (e_\omega \cdot \Theta)}$, for an ideal
${\mathfrak A}'_M$ of $M$ since 
$\Theta = \frac{\sigma}{2}(1-\sigma)$; then using the relation
$e_\omega \cdot \Theta = \Theta \cdot e_{\chi_0}$, we obtain
$(y^{e_\omega}) = {\mathfrak A}_L^\Theta
={\mathfrak A}'^{1-\sigma}_L$ for an ideal ${\mathfrak A}'_L$ of $L$. 

\smallskip
We know that from the relation
$\zeta_3 = {\rm N}_{M/K}(y^{e_\omega})$, with
$(y^{e_\omega}) = {\mathfrak A}'^{1-\sigma}_L$, we must consider
the ideal $(\alpha) = {\rm N}_{M/K}({\mathfrak A}'_L)$ of $K$, which
gives rise to 
$\Lambda_1 = \langle \zeta_3, \ell_1, \ldots, \ell_n, \alpha \rangle$.
But $\alpha \in \Q^\times$ yielding 
$(\Lambda_1 : \Lambda_1 \cap {\mathcal N}_{M/K})=1$.
So $\order (\Cl_2/\Cl_1) = 3^{n-d_3} =  3^{n-1} \equiv 0 \pmod 3$ 
(absurd). Thus $\Cl_M = 1$.
\end{proof}

\begin{remarks} 
(i) When $\Cl_M=1$, the relation $n = \nu_3+d_3$ gives 
easily the solutions $N$.

\smallskip
\quad (ii) If $N=\prod_{i=1}^n \ell_i$, with $\ell_i \equiv 1 \pmod 3$, 
$\nu_3=1$ ($\zeta_3$ non-norm)
and $N \equiv 1 \pmod 9$ ($3$ unramified), one finds easily that 
$\Cl_2/\Cl_1\ne 1$ giving ${\rm rk}_3(\Cl_M) \geq 2$ (e.g., 
$N=13\cdot 97, 7 \cdot 283, 7 \cdot 337$ for which 
$\Cl_M \simeq (\Z/3\Z)^2$). For similar results using the
properties of the filtration with Hilbert's symbols, see 
\cite{AAIMT,ATIA1,ATIA2,Gerth,Gras0,Gras8,Kobayashi}.

\smallskip
\quad (iii) Using for any $p$ the fixed point formula for ${\mathcal T}_M$, 
it is possible to compute the sets of primes $\ell$ dividing $N$ for which the 
set $\{{\mathfrak l}, \ {\mathfrak l} \mid \ell\}$ is $p$-primitif or to use
Theorem \ref{equiv}.
But the equivalence: ``${\mathcal T}_L$ is $p$-rational if and only if 
${\mathcal T}_M$ is $p$-rational'' does not work for $p>3$ as
shown by the following example:

\smallskip
\footnotesize
\begin{verbatim}
N=11=[11,1]  5 ramified      
Structure of C_L=[5]          Structure of C_M=[5,5]
Structure of T_L=[]           Structure of T_M=[5]
\end{verbatim}

\smallskip
\normalsize
Same remark for the trivialities of $\Cl_M$ and $\Cl_L$:

\smallskip
\footnotesize
\begin{verbatim}
N=357=[3,1;7,1;17,1]  5 unramified
q=3 ClassqL=[3]       q=3 ClassqM=[3,3,3,3]       q=5 ClassqM=[5]
\end{verbatim}
\normalsize
\end{remarks}

\normalsize
\subsection{The logarithmic class group}\label{bp}
The logarithmic class group $\wt \Cl_F$ of
a number field $F$ was introduced by Jaulent in
\cite{Jaulent0,Jalog} in the following context.

\subsubsection{The fundamental exact sequence}
There exists an exact sequence, defining the kernel $\wt \Cl_F^{[p]}$, of the form:
$$1\too \wt \Cl_F^{[p]} \tooo \wt \Cl_F \tooo \Cl_F^{P_F}  \too 1, $$

where $\Cl_F^{P_F}  := \Cl_F/ \langle \cl_F({\mathfrak p}),\ 
{\mathfrak p} \mid p \rangle$.
The logarithmic class group is given by the instruction
${\sf bnflog(K,p)}$; warning, the result is written in the following order
(see \cite{BJ}): $[ \wt \Cl_F],\  [\wt \Cl_F^{[p]}],\  [\Cl_F^{P_F} ]$.

\smallskip
For the field $M$, if $p$ is totally ramified, then $\Cl_M^{P_M}  := 
\Cl_M / \langle \cl_M({\mathfrak P}_M) \rangle$ where 
$(p)={\mathfrak P}_M^p$;
thus $0 \leq {\rm rk}_p(\Cl_M) - {\rm rk}_p(\Cl_M^{P_M}) \leq 1$;
otherwise, if $p$ splits in $M/K$, one gets
$0 \leq {\rm rk}_p(\Cl_M) - {\rm rk}_p(\Cl_M^{P_M}) \leq p-1$.

\smallskip
For instance we will obtain:

\footnotesize
\smallskip
\begin{verbatim}
N=101=[101,1] 5 unramified          N=124=[2,2;31,1] 5 unramified
Structure of C_L=[5]                Structure of C_L=[5] 
Structure of T_L=[25]               Structure of T_L=[125] 
Clog_L=[],[],[]                     Clog_L=[25],[5],[5]
Structure of C_M=[5,5,5,5]          Structure of C_M=[5,5] 
Structure of T_M=[25,25,5,5,5]      Structure of T_M=[125,25,25,25,5,5]
Clog_M=[5,5,5],[5,5],[5]            Clog_M=[25,25,5,5,5],[25,5,5,5],[5,5]
\end{verbatim}

\smallskip
\normalsize
In other words, for $N=124$, we have $\wt \Cl_M= \!\sf [25,25,5,5,5]$,
$\wt \Cl_M^{[p]}= \sf [25,5,5,5]$, $\Cl_M^{P_M} = \sf [5,5]$,
considering that $\Cl_M = \sf [5,5]$. 

\subsubsection{Reflection theorem for the logarithmic class group}
From \cite[Scolie 3.13]{Jaulent0}, when $F$ contains $\mu_p$ 
and for an automorphism group $g$ of $F$, of prime to $p$ order, 
we have, for all $\chi$ (under the Leopoldt conjecture):
\begin{equation}\label{eq7}
{\rm rk}_{\chi}(\wt \Cl_F) = {\rm rk}_{\chi*}({\mathcal T}_F^{\rm bp}),
\end{equation}

where ${\mathcal T}_F^{\rm bp}$ is the Bertrandias-Payan module
(see Lemma \ref{exact}). We then have (cf. Schema of Lemma \ref{exact}): 
$${\rm rk}_{\chi}(\wt \Cl_F) = 
{\rm rk}_{\chi^*}({\mathcal T}_F^{\rm bp}) \leq
{\rm rk}_{\chi^*}( \ov \Cl_F) +  {\rm rk}_{\chi^*}( {\mathcal R}_F)$$
and:
$${\rm rk}_{\chi}(\wt \Cl_F) \leq
{\rm rk}_{\chi}(\wt \Cl_F^{[p]}) +  {\rm rk}_{\chi}(\Cl_F/\cl_F(P_F)). $$

\smallskip
The relation \eqref{eq7} applies in $M/L$ and may give some information:
\begin{equation}
\begin{aligned}
{\rm rk}_{\chi}(\wt \Cl_M) &= {\rm rk}_{\chi*}({\mathcal T}_M^{\rm bp}),\\
{\rm rk}_{\omega}(\wt \Cl_M)& = {\rm rk}_p({\mathcal T}_L^{\rm bp}),\\
{\rm rk}_p(\wt \Cl_L) &= {\rm rk}_{\omega}({\mathcal T}_M^{\rm bp}),\\
{\rm rk}_p(\wt \Cl_M)& = {\rm rk}_p({\mathcal T}_M^{\rm bp}).
\end{aligned}
\end{equation}

For instance in the most usual case $p$ totally ramified ($d_p=0$), we know that
${\rm rk}_p({\mathcal W}_M)={\rm rk}_p({\mathcal W}_L)=0$, which yields
${\mathcal T}_M^{\rm bp}={\mathcal T}_M$ then (Corollary \ref{ranks}):
\begin{equation}
\begin{aligned}
{\rm rk}_{\chi}(\wt \Cl_M) &= {\rm rk}_{\chi*}({\mathcal T}_M) \! =
{\rm rk}_\chi(\Cl_M/\cl_M(P_M)) , \\
{\rm rk}_{\omega}(\wt \Cl_M)& = {\rm rk}_p({\mathcal T}_L) \ = 
{\rm rk}_\omega(\Cl_M/\cl_M(P_M)), \\
{\rm rk}_p(\wt \Cl_L) &= {\rm rk}_{\omega}({\mathcal T}_M) =
{\rm rk}_p(\Cl_L/\cl_L(P_L)),\\
{\rm rk}_p(\wt \Cl_M)& = {\rm rk}_p({\mathcal T}_M) =
{\rm rk}_p(\Cl_M/\cl_M(P_M)).
\end{aligned}
\end{equation}

\smallskip
When $p$ splits in $M/K$ one obtains for instance, 
since ${\rm rk}_{\chi}({\mathcal W}_M) =1$:
\begin{equation*}
{\rm rk}_p(\wt \Cl_L)  \geq {\rm rk}_{\omega}({\mathcal T}_M) 
- {\rm rk}_{\omega}({\mathcal W}_M) = 
{\rm rk}_{\omega}({\mathcal T}_M) - 1.
\end{equation*}

One may examine the above relations with the following example:

\footnotesize
\smallskip
\begin{verbatim}
N=182=[2,1;7,1;13,1]  3 ramified          N=99=[3,2;11,1]  5 unramified
Structure of C_L=[3,3,3]                  Structure of C_L=[5]
Structure of T_L=[3,3]                    Structure of T_L=[5]
Clog_L=[3,3],[],[3,3]                     Clog_L=[],[],[]
Structure of C_M=[3,3,3,3,3]              Structure of C_M=[5,5]
Structure of T_M=[3,3,3,3]                Structure of T_M=[5,5,5,5,5]
Clog_M=[3,3,3,3],[],[3,3,3,3]             Clog_M=[5],[],[5]
\end{verbatim}

\medskip
\normalsize
All these relations are confirmed by the forthcomming tables.

\newpage
\appendix 

\section{Programs and Numerical results}\label{A}
We give first the general PARI/GP \cite{pari} program computing
the $p$-structure of the groups $\Cl_L$, ${\mathcal T}_L$, $\wt \Cl_L$,
$\Cl_M$, ${\mathcal T}_M$, $\wt \Cl_M$, then various tables.
Sometimes, the groups $\wt \Cl$ are not computable by PARI/GP
in a reasonable running time.

\subsection{PARI/GP program}\label{T}
Since $\Q(\sqrt[p]{N}) = \Q(\sqrt[p]{N'})$ if and only if
$N' = N^a A^p$, $a \in \Z$ prime to $p$, $A \in \Q^\times$,
the program computes the minimal representative 
for each field.
The fields $L$, $K$ are defined with the polynomials
$x^p-N$ and the $p$th cyclotomic one, respectively; the 
field $M$ is defined by the function ${\sf polcompositum}$
of PARI/GP.
The computation of ${\mathcal T}_L$ and ${\mathcal T}_M$
use ray class fields with modulus $p^n$ for $n$ large enough
such that ${\rm Gal}(H_L^{(p^n)}/L)=: 
(p^{a_1},\cdots ,p^{a_{(p-1)/2+1}} ;  p^{a_1},\cdots ,p^{a_{R_L}})$
where the $(p-1)/2+1$ large components are relative to the ``free part'',
$( p^{a_1},\cdots ,p^{a_{R_L}})$ giving the $p$-structure of 
${\mathcal T}_L$. So one must verify that  $n$ is much larger than the
exponent of ${\mathcal T}_L$. For the computation of ${\mathcal T}_M$,
replace $(p-1)/2+1$ by $p\,(p-1)/2+1$.

\medskip
\footnotesize
\begin{verbatim}
{p=3;n=8;for(N=2,10^4,F=factor(N);F1=component(F,1);F2=component(F,2);
d=matsize(F1)[1]; T=1;for(j=1,d,e=F2[j];if(e>=p,T=0));if(T==0,next); 
T=1;for(k=2,p-1,NN=1;for(j=1,d,ek=lift(Mod(k*F2[j],p)); 
NN=NN*F1[j]^ek);if(NN<N,T=0;break));if(T==0,next);print();
if(Mod(N,p)!=0 & Mod(N^(p-1),p^2)==1,
print("N=",N,"=",factor(N)," ",p," unramified"));
if(Mod(N,p)==0 || Mod(N^(p-1),p^2)!=1,
print("N=",N,"=",factor(N)," ",p," ramified"));
PL=x^p-N;L=bnfinit(PL,1);PK=polcyclo(p); 
PM=polcompositum(PL,PK)[1];M=bnfinit(PM,1); 
ClassL=L.cyc; \\ print(ClassL); 
ClogL=bnflog(L,p); \\ logarithmic class group of L
L0=List;e=matsize(ClassL)[2];for(k=1,e,c=ClassL[e-k+1]; 
if(Mod(c,p)==0,listinsert(L0,p^valuation(c,p),1))); 
ClassL=L0; \\ structure of the p-class group of L 
ClassM=M.cyc; \\ print(ClassM); 
M0=List;e=matsize(ClassM)[2];for(k=1,e,c=ClassM[e-k+1]; 
if(Mod(c,p)==0,listinsert(M0,p^valuation(c,p),1))); 
ClassM=M0; \\ structure of the p-class group of M 
ClogM=bnflog(M,p); \\ logarithmic class group of M
Lpn=bnrinit(L,p^n);HLpn=Lpn.cyc; \\ structure of the ray class field over L 
TorL=List;e=matsize(HLpn)[2];for(k=1,e-(p-1)/2-1,c=HLpn[e-k+1]; 
if(Mod(c,p)==0,listinsert(TorL,p^valuation(c,p),1)));
print("Structure of C_L=",ClassL);print("Structure of T_L=",TorL); 
print("Structure of Clog_L=",ClogL);
Mpn=bnrinit(M,p^n);HMpn=Mpn.cyc; \\ structure of the ray class field over M 
TorM=List;e=matsize(HMpn)[2];for(k=1,e-p*(p-1)/2-1,c=HMpn[e-k+1]; 
if(Mod(c,p)==0,listinsert(TorM,p^valuation(c,p),1)));
print("Structure of C_M=",ClassM);print("Structure of T_M=",TorM);
print("Structure of Clog_M=",ClogM))}
\end{verbatim}
\normalsize

\medskip
Since the program computes the whole structure of the class 
groups of $L$ and $M$, the reader may print these results replacing:

\begin{verbatim}
\\ print(ClassL); \\ print(ClassM);
\end{verbatim}
by:
\begin{verbatim}
print(ClassL);print(ClassM);
\end{verbatim}

\textwidth=35pc
\subsection{Table of $\Cl_L$, ${\mathcal T}_L$, $\wt \Cl_L$,
$\Cl_M$, ${\mathcal T}_M$, $\wt \Cl_M$, $p = 3$}
${}$
\scriptsize
\textheight=254mm
\twocolumn

\normalsize
\onecolumn
\scriptsize
\begin{verbatim}
N=9996=[2,2;3,1;7,2;17,1] 3 ramified
Structure of C_L=[3,3]
Structure of T_L=[3,3]
Clog_L=[3],[],[3]
Structure of C_M=[3,3,3,3]
Structure of T_M=[3,3,3]
Clog_M=[3,3,3],[],[3,3,3]

N=9997=[13,1;769,1] 3 ramified
Structure of C_L=[3,3]
Structure of T_L=[3,3]
Clog_L=[3],[],[3]
Structure of C_M=[3,3,3,3]
Structure of T_M=[3,3,3]
Clog_M=[3,3,3],[],[3,3,3]

N=9998=[2,1;4999,1] 3 unramified
Structure of C_L=[3]
Structure of T_L=[3]
Clog_L=[],[],[]
Structure of C_M=[3]
Structure of T_M=[3,3]
Clog_M=[],[],[]

N=9999=[3,2;11,1;101,1] 3 ramified
Structure of C_L=[3]
Structure of T_L=[3]
Clog_L=[],[],[]
Structure of C_M=[3,3]
Structure of T_M=[3]
Clog_M=[3],[],[3]
\end{verbatim}
\normalsize

\newpage
\subsection{Table of $\Cl_L$, ${\mathcal T}_L$, $\wt \Cl_L$,
$\Cl_M$, ${\mathcal T}_M$, $\wt \Cl_M$, $p=5$}
${}$
\scriptsize
\medskip
\begin{verbatim}
N=2=[2,1] 5 ramified
Structure of C_L=[]
Structure of T_L=[]
Clog_L=[],[],[]
Structure of C_M=[]
Structure of T_M=[]
Clog_M=[],[],[]

N=3=[3,1] 5 ramified
Structure of C_L=[]
Structure of T_L=[]
Clog_L=[],[],[]
Structure of C_M=[]
Structure of T_M=[]
Clog_M=[],[],[]

N=5=[5,1] 5 ramified
Structure of C_L=[]
Structure of T_L=[]
Clog_L=[],[],[]
Structure of C_M=[]
Structure of T_M=[]
Clog_M=[],[],[]

N=6=[2,1;3,1] 5 ramified
Structure of C_L=[]
Structure of T_L=[5]
Clog_L=[],[],[]
Structure of C_M=[5]
Structure of T_M=[5]
Clog_M=[5],[],[5]

N=7=[7,1] 5 unramified
Structure of C_L=[]
Structure of T_L=[25]
Clog_L=[],[],[]
Structure of C_M=[]
Structure of T_M=[25,5,5,5]
Clog_M=[5],[5],[]

N=10=[2,1;5,1] 5 ramified
Structure of C_L=[]
Structure of T_L=[]
Clog_L=[],[],[]
Structure of C_M=[]
Structure of T_M=[]
Clog_M=[],[],[]

N=11=[11,1] 5 ramified
Structure of C_L=[5]
Structure of T_L=[]
Clog_L=[],[],[]
Structure of C_M=[5,5]
Structure of T_M=[5]
Clog_M=[5],[],[5]

N=12=[2,2;3,1] 5 ramified
Structure of C_L=[]
Structure of T_L=[5]
Clog_L=[],[],[]
Structure of C_M=[5]
Structure of T_M=[5]
Clog_M=[5],[],[5]
\end{verbatim}
\normalsize
\twocolumn
\scriptsize

\normalsize

\onecolumn
\subsection{Table of $\Cl_L$, ${\mathcal T}_L$, $\wt \Cl_L$,
$\Cl_M$, ${\mathcal T}_M$, $\wt \Cl_M$, $p=5$, 
$N \equiv 1 \!\!\pmod 5$}
${}$
\scriptsize
\begin{verbatim}
N=11=[11,1]
Structure of C_L=[5]
Structure of T_L =[]
Structure of C_M=[5,5]
Structure of T_M=[5]

N=31=[31,1]
Structure of C_L=[5,5]
Structure of T_L =[5]
Structure of C_M=[5,5,5,5,5]
Structure of T_M=[5,5,5,5]

N=41=[41,1]
Structure of C_L=[5]
Structure of T_L =[]
Structure of C_M=[5,5]
Structure of T_M=[5]

N=61=[61,1]
Structure of C_L=[5]
Structure of T_L =[]
Structure of C_M=[5,5]
Structure of T_M=[5]

N=71=[71,1]
Structure of C_L=[5]
Structure of T_L =[]
Structure of C_M=[5,5]
Structure of T_M=[5]

N=101=[101,1]
Structure of C_L=[5]
Structure of T_L =[25]
Structure of C_M=[5,5,5,5]
Structure of T_M=[25,25,5,5,5]

N=121=[11,2]
Structure of C_L=[5]
Structure of T_L =[]
Structure of C_M=[5,5]
Structure of T_M=[5]

N=131=[131,1]
Structure of C_L=[5,5]
Structure of T_L =[5]
Structure of C_M=[5,5,5,5]
Structure of T_M=[5,5,5]

N=151=[151,1]
Structure of C_L=[5]
Structure of T_L =[25]
Structure of C_M=[5,5]
Structure of T_M=[25,5,5,5,5]

N=181=[181,1]
Structure of C_L=[5,5]
Structure of T_L =[5]
Structure of C_M=[5,5,5,5]
Structure of T_M=[5,5,5]

N=191=[191,1]
Structure of C_L=[5]
Structure of T_L =[5]
Structure of C_M=[5,5,5,5]
Structure of T_M=[5,5,5,5]
\end{verbatim}
\normalsize
\twocolumn
\scriptsize
\begin{verbatim}
N=211=[211,1]
Structure of C_L=[5,5,5]
Structure of T_L =[5,5]
Structure of C_M=[5,5,5,5,5,5,5,5,5]
Structure of T_M=[5,5,5,5,5,5,5,5]

N=241=[241,1]
Structure of C_L=[5]
Structure of T_L =[]
Structure of C_M=[5,5]
Structure of T_M=[5]

N=251=[251,1]
Structure of C_L=[5]
Structure of T_L =[25]
Structure of C_M=[5,5]
Structure of T_M=[25,5,5,5,5,5]

N=271=[271,1]
Structure of C_L=[5]
Structure of T_L =[5]
Structure of C_M=[5,5,5,5]
Structure of T_M=[5,5,5,5]

N=281=[281,1]
Structure of C_L=[5,5,5]
Structure of T_L =[5,5]
Structure of C_M=[25,5,5,5,5,5,5,5]
Structure of T_M=[5,5,5,5,5,5,5]

N=311=[311,1]
Structure of C_L=[5]
Structure of T_L =[]
Structure of C_M=[5,5]
Structure of T_M=[5]

N=331=[331,1]
Structure of C_L=[5]
Structure of T_L =[]
Structure of C_M=[5,5]
Structure of T_M=[5]

N=341=[11,1;31,1]
Structure of C_L=[5,5,5]
Structure of T_L =[5,5]
Structure of C_M=[5,5,5,5,5,5,5,5,5]
Structure of T_M=[5,5,5,5,5,5,5,5]

N=401=[401,1]
Structure of C_L=[5,5]
Structure of T_L =[25]
Structure of C_M=[5,5,5,5,5]
Structure of T_M=[25,25,5,5,5,5]

N=421=[421,1]
Structure of C_L=[5]
Structure of T_L =[]
Structure of C_M=[5,5,5]
Structure of T_M=[5,5]

N=431=[431,1]
Structure of C_L=[5]
Structure of T_L =[]
Structure of C_M=[5,5]
Structure of T_M=[5]

N=451=[11,1;41,1]
Structure of C_L=[5,5]
Structure of T_L =[25]
Structure of C_M=[5,5,5,5,5,5]
Structure of T_M=[25,5,5,5,5,5,5,5]

N=461=[461,1]
Structure of C_L=[5]
Structure of T_L =[]
Structure of C_M=[5,5,5]
Structure of T_M=[5,5]

N=491=[491,1]
Structure of C_L=[5]
Structure of T_L =[]
Structure of C_M=[5,5]
Structure of T_M=[5]

N=521=[521,1]
Structure of C_L=[5]
Structure of T_L =[]
Structure of C_M=[5,5,5]
Structure of T_M=[5,5]

N=541=[541,1]
Structure of C_L=[5,5]
Structure of T_L =[5]
Structure of C_M=[5,5,5,5]
Structure of T_M=[5,5,5]

N=571=[571,1]
Structure of C_L=[5,5]
Structure of T_L =[5]
Structure of C_M=[5,5,5,5]
Structure of T_M=[5,5,5]

N=601=[601,1]
Structure of C_L=[5]
Structure of T_L =[25]
Structure of C_M=[5,5]
Structure of T_M=[25,5,5,5,5,5]

N=631=[631,1]
Structure of C_L=[5]
Structure of T_L =[]
Structure of C_M=[5,5]
Structure of T_M=[5]

N=641=[641,1]
Structure of C_L=[5]
Structure of T_L =[5]
Structure of C_M=[5,5,5,5]
Structure of T_M=[5,5,5,5]

N=661=[661,1]
Structure of C_L=[5]
Structure of T_L =[]
Structure of C_M=[5,5]
Structure of T_M=[5]

N=671=[11,1;61,1]
Structure of C_L=[5,5,5]
Structure of T_L =[5,5]
Structure of C_M=[5,5,5,5,5,5,5,5]
Structure of T_M=[5,5,5,5,5,5,5]

N=691=[691,1]
Structure of C_L=[5,5,5]
Structure of T_L =[5,5]
Structure of C_M=[5,5,5,5,5,5,5,5]
Structure of T_M=[25,5,5,5,5,5,5]

N=701=[701,1]
Structure of C_L=[5,5]
Structure of T_L =[25]
Structure of C_M=[5,5,5,5,5]
Structure of T_M=[25,25,5,5,5,5]

N=751=[751,1]
Structure of C_L=[5,5]
Structure of T_L =[25,5]
Structure of C_M=[5,5,5,5]
Structure of T_M=[25,25,25,5,5,5,5]

N=761=[761,1]
Structure of C_L=[5,5]
Structure of T_L =[5]
Structure of C_M=[5,5,5,5,5]
Structure of T_M=[5,5,5,5]

N=781=[11,1;71,1]
Structure of C_L=[5,5,5]
Structure of T_L =[5,5]
Structure of C_M=[5,5,5,5,5,5,5,5,5]
Structure of T_M=[5,5,5,5,5,5,5,5]

N=811=[811,1]
Structure of C_L=[5]
Structure of T_L =[]
Structure of C_M=[5,5]
Structure of T_M=[5]

N=821=[821,1]
Structure of C_L=[5]
Structure of T_L =[]
Structure of C_M=[5,5]
Structure of T_M=[5]

N=881=[881,1]
Structure of C_L=[5]
Structure of T_L =[]
Structure of C_M=[5,5,5]
Structure of T_M=[5,5]

N=911=[911,1]
Structure of C_L=[5]
Structure of T_L =[]
Structure of C_M=[5,5]
Structure of T_M=[5]

N=941=[941,1]
Structure of C_L=[5,5]
Structure of T_L =[5]
Structure of C_M=[5,5,5,5]
Structure of T_M=[5,5,5]

N=961=[31,2]
Structure of C_L=[5,5]
Structure of T_L =[5]
Structure of C_M=[5,5,5,5,5]
Structure of T_M=[5,5,5,5]

N=971=[971,1]
Structure of C_L=[5,5]
Structure of T_L =[5]
Structure of C_M=[5,5,5,5]
Structure of T_M=[5,5,5]

N=991=[991,1]
Structure of C_L=[5]
Structure of T_L =[]
Structure of C_M=[5,5,5]
Structure of T_M=[5,5]

N=1021=[1021,1]
Structure of C_L=[5]
Structure of T_L =[]
Structure of C_M=[5,5]
Structure of T_M=[5]

N=1031=[1031,1]
Structure of C_L=[5]
Structure of T_L =[]
Structure of C_M=[5,5,5]
Structure of T_M=[5,5]

N=1051=[1051,1]
Structure of C_L=[5]
Structure of T_L =[25]
Structure of C_M=[5,5]
Structure of T_M=[25,5,5,5,5]

N=1061=[1061,1]
Structure of C_L=[5]
Structure of T_L =[]
Structure of C_M=[5,5]
Structure of T_M=[5]
\end{verbatim}
\normalsize

\onecolumn
\subsection{Table of $\Cl_L$, ${\mathcal T}_L$, $\Cl_R$,
${\mathcal T}_R$, $\Cl_M$, ${\mathcal T}_M$,  $p=5$, $R=\Q(\sqrt 5)$}
${}$
\textheight=234mm
\scriptsize
\begin{verbatim}
N=2=[2,1] 5 ramified
Structure of C_L=[]
Structure of T_L=[]
Structure of C_R=[]
Structure of T_R=[]
Structure of C_M=[]
Structure of T_M=[]

N=3=[3,1] 5 ramified
Structure of C_L=[]
Structure of T_L=[]
Structure of C_R=[]
Structure of T_R=[]
Structure of C_M=[]
Structure of T_M=[]

N=5=[5,1] 5 ramified
Structure of C_L=[]
Structure of T_L=[]
Structure of C_R=[]
Structure of T_R=[]
Structure of C_M=[]
Structure of T_M=[]

N=6=[2,1;3,1] 5 ramified
Structure of C_L=[]
Structure of T_L=[5]
Structure of C_R=[]
Structure of T_R=[]
Structure of C_M=[5]
Structure of T_M=[5]

N=7=[7,1] 5 unramified
Structure of C_L=[]
Structure of T_L=[25]
Structure of C_R=[]
Structure of T_R=[]
Structure of C_M=[]
Structure of T_M=[25,5,5,5]

N=10=[2,1;5,1] 5 ramified
Structure of C_L=[]
Structure of T_L=[]
Structure of C_R=[]
Structure of T_R=[]
Structure of C_M=[]
Structure of T_M=[]

N=11=[11,1] 5 ramified
Structure of C_L=[5]
Structure of T_L=[]
Structure of C_R=[5]
Structure of T_R=[]
Structure of C_M=[5,5]
Structure of T_M=[5]

N=12=[2,2;3,1] 5 ramified
Structure of C_L=[]
Structure of T_L=[5]
Structure of C_R=[]
Structure of T_R=[]
Structure of C_M=[5]
Structure of T_M=[5]

N=13=[13,1] 5 ramified
Structure of C_L=[]
Structure of T_L=[]
Structure of C_R=[]
Structure of T_R=[]
Structure of C_M=[]
Structure of T_M=[]
\end{verbatim}
\textheight=254mm
\twocolumn

\normalsize

\onecolumn
\subsection{Table of $\Cl_L$, ${\mathcal T}_L$, $\wt \Cl_L$,
$\Cl_M$, ${\mathcal T}_M$, $\wt \Cl_M$, $p=7$}
${}$
\textheight=234mm
\scriptsize
\begin{verbatim}
N=2=[2,1] 7 ramified
Structure of C_L=[]
Structure of T_L=[]
Clog_L=[[],[],[]]
Structure of C_M=[]
Structure of T_M=[]
Clog_M=[[],[],[]]

N=3=[3,1] 7 ramified
Structure of C_L=[]
Structure of T_L=[]
Clog_L=[[],[],[]]
Structure of C_M=[]
Structure of T_M=[]
Clog_M=[[],[],[]]

N=5=[5,1] 7 ramified
Structure of C_L=[]
Structure of T_L=[]
Clog_L=[[],[],[]]
Structure of C_M=[]
Structure of T_M=[]
Clog_M=[[],[],[]]

N=6=[2,1;3,1] 7 ramified
Structure of C_L=[]
Structure of T_L=[7]
Clog_L=[[],[],[]]
Structure of C_M=[7]
Structure of T_M=[7]
Clog_M=[[7],[],[7]]

N=7=[7,1] 7 ramified
Structure of C_L=[]
Structure of T_L=[]
Clog_L=[[],[],[]]
Structure of C_M=[]
Structure of T_M=[]
Clog_M=[[],[],[]]

N=10=[2,1;5,1] 7 ramified
Structure of C_L=[]
Structure of T_L=[7]
Clog_L=[[],[],[]]
Structure of C_M=[7]
Structure of T_M=[7]
Clog_M=[[7],[],[7]]

N=11=[11,1] 7 ramified
Structure of C_L=[]
Structure of T_L=[]
Clog_L=[[],[],[]]
Structure of C_M=[]
Structure of T_M=[]
Clog_M=[[],[],[]]

N=12=[2,2;3,1] 7 ramified
Structure of C_L=[]
Structure of T_L=[7]
Clog_L=[[],[],[]]
Structure of C_M=[7]
Structure of T_M=[7]
Clog_M=[[7],[],[7]]
\end{verbatim}
\textheight=254mm
\twocolumn
\begin{verbatim}
N=13=[13,1] 7 ramified
Structure of C_L=[7]
Structure of T_L=[]
Clog_L=[[],[],[]]
Structure of C_M=[7,7,7]
Structure of T_M=[7,7]
Clog_M=[[7,7],[],[7,7]]

N=14=[2,1;7,1] 7 ramified
Structure of C_L=[]
Structure of T_L=[]
Clog_L=[[],[],[]]
Structure of C_M=[]
Structure of T_M=[]
Clog_M=[[],[],[]]

N=15=[3,1;5,1] 7 ramified
Structure of C_L=[]
Structure of T_L=[7]
Clog_L=[[],[],[]]
Structure of C_M=[7]
Structure of T_M=[7]
Clog_M=[[7],[],[7]]

N=17=[17,1] 7 ramified
Structure of C_L=[]
Structure of T_L=[]
Clog_L=[[],[],[]]
Structure of C_M=[]
Structure of T_M=[]
Clog_M=[[],[],[]]

N=18=[2,1;3,2] 7 unramified
Structure of C_L=[]
Structure of T_L=[7]
Clog_L=[[],[],[]]
Structure of C_M=[]
Structure of T_M=[7,7,7,7,7,7]
Clog_M=[[],[],[]]

N=19=[19,1] 7 unramified
Structure of C_L=[]
Structure of T_L=[49]
Clog_L=[[],[],[]]
Structure of C_M=[]
Structure of T_M=[49,49,49,7,7,7]
Clog_M=[[7,7,7],[7,7,7],[]]

N=20=[2,2;5,1] 7 ramified
Structure of C_L=[]
Structure of T_L=[7]
Clog_L=[[],[],[]]
Structure of C_M=[7]
Structure of T_M=[7]
Clog_M=[[7],[],[7]]

N=21=[3,1;7,1] 7 ramified
Structure of C_L=[]
Structure of T_L=[]
Clog_L=[[],[],[]]
Structure of C_M=[]
Structure of T_M=[]
Clog_M=[[],[],[]]

N=22=[2,1;11,1] 7 ramified
Structure of C_L=[7]
Structure of T_L=[7]
Clog_L=[[],[],[]]
Structure of C_M=[7,7,7,7]
Structure of T_M=[7,7,7]
Clog_M=[[7,7,7],[],[7,7,7]]

N=23=[23,1] 7 ramified
Structure of C_L=[]
Structure of T_L=[]
Clog_L=[[],[],[]]
Structure of C_M=[]
Structure of T_M=[]
Clog_M=[[],[],[]]

N=24=[2,3;3,1] 7 ramified
Structure of C_L=[]
Structure of T_L=[7]
Clog_L=[[],[],[]]
Structure of C_M=[7]
Structure of T_M=[7]
Clog_M=[[7],[],[7]]

N=26=[2,1;13,1] 7 ramified
Structure of C_L=[7]
Structure of T_L=[7]
Clog_L=[[],[],[]]
Structure of C_M=[7,7,7,7]
Structure of T_M=[7,7,7]
Clog_M=[[7,7,7],[],[7,7,7]]

N=28=[2,2;7,1] 7 ramified
Structure of C_L=[]
Structure of T_L=[]
Clog_L=[[],[],[]]
Structure of C_M=[]
Structure of T_M=[]
Clog_M=[[],[],[]]

N=29=[29,1] 7 ramified
Structure of C_L=[7]
Structure of T_L=[]
Clog_L=[[],[],[]]
Structure of C_M=[7,7,7]
Structure of T_M=[7,7]
Clog_M=[[7,7],[],[7,7]]

N=30=[2,1;3,1;5,1] 7 unramified
Structure of C_L=[]
Structure of T_L=[7,7]
Clog_L=[[],[],[]]
Structure of C_M=[7]
Structure of T_M=[7,7,7,7,7,7,7]
Clog_M=[[7],[],[7]]

N=31=[31,1] 7 unramified
Structure of C_L=[]
Structure of T_L=[49]
Clog_L=[[],[],[]]
Structure of C_M=[]
Structure of T_M=[49,7,7,7,7,7]
Clog_M=[[7],[7],[]]

N=33=[3,1;11,1] 7 ramified
Structure of C_L=[]
Structure of T_L=[7]
Clog_L=[[],[],[]]
Structure of C_M=[7]
Structure of T_M=[7]
Clog_M=[[7],[],[7]]
\end{verbatim}
\normalsize

\onecolumn
\section{Tables of $\Cl_L$, ${\mathcal T}_L$, 
$p=7, 11, 13$}\label{B}
We only write the non-trivial structures of $\Cl_L$, ${\mathcal T}_L$
in the case $N$ prime, $N \equiv 1 \pmod p$:
${}$
\smallskip
\scriptsize
\begin{verbatim}
{p=7;n=4;forprime(N=1,10^4,if(Mod(N,p)!=1,next);
PL=x^p-N;L=bnfinit(PL,1);ClassL=L.cyc;
L0=List;e=matsize(ClassL)[2];for(k=1,e,c=ClassL[e-k+1];
if(Mod(c,p)==0,listinsert(L0,p^valuation(c,p),1)));
ClassL=L0;\\structureofthep-classgroup
Lpn=bnrinit(L,p^n);\\rayclassfieldmodp^n
HLpn=Lpn.cyc;\\structureoftherayclassfield
TorL=List;e=matsize(HLpn)[2];for(k=1,e-(p-1)/2-1,c=HLpn[e-k+1];
if(Mod(c,p)==0,listinsert(TorL,p^valuation(c,p),1)));
if(ClassL!=[] || TorL!=[],print();print("N=",N,"=",factor(N));
print("Structure of C_L=",ClassL);print("Structure of T_L=",TorL)))}
\end{verbatim}
\normalsize

\smallskip
\subsection{Table of $\Cl_L$, ${\mathcal T}_L$, $p =7$}
${}$
\scriptsize
\begin{verbatim}
N=29=[29,1]
Structure of C_L=[7]
Structure of T_L=[]

N=43=[43,1]
Structure of C_L=[7]
Structure of T_L=[7]

N=71=[71,1]
Structure of C_L=[7]
Structure of T_L=[]

N=113=[113,1]
Structure of C_L=[7]
Structure of T_L=[]

N=127=[127,1]
Structure of C_L=[7,7]
Structure of T_L=[7]

N=197=[197,1]
Structure of C_L=[7]
Structure of T_L=[49]

N=211=[211,1]
Structure of C_L=[7]
Structure of T_L=[]

N=239=[239,1]
Structure of C_L=[7]
Structure of T_L=[]

N=281=[281,1]
Structure of C_L=[7]
Structure of T_L=[7]

N=337=[337,1]
Structure of C_L=[7,7]
Structure of T_L=[7]

N=379=[379,1]
Structure of C_L=[7]
Structure of T_L=[]

N=421=[421,1]
Structure of C_L=[7]
Structure of T_L=[]
\end{verbatim}
\normalsize
\textheight=242mm
\twocolumn
\scriptsize
\begin{verbatim}
N=449=[449,1]
Structure of C_L=[7]
Structure of T_L=[]

N=463=[463,1]
Structure of C_L=[7]
Structure of T_L=[]

N=491=[491,1]
Structure of C_L=[7]
Structure of T_L=[49]

N=547=[547,1]
Structure of C_L=[7]
Structure of T_L=[]

N=617=[617,1]
Structure of C_L=[7]
Structure of T_L=[]

N=631=[631,1]
Structure of C_L=[7,7]
Structure of T_L=[7]

N=659=[659,1]
Structure of C_L=[7,7]
Structure of T_L=[7]

N=673=[673,1]
Structure of C_L=[7]
Structure of T_L=[]

N=701=[701,1]
Structure of C_L=[7]
Structure of T_L=[]

N=743=[743,1]
Structure of C_L=[7]
Structure of T_L=[]

N=757=[757,1]
Structure of C_L=[7]
Structure of T_L=[]

N=827=[827,1]
Structure of C_L=[7]
Structure of T_L=[]

N=883=[883,1]
Structure of C_L=[7]
Structure of T_L=[49]

N=911=[911,1]
Structure of C_L=[7,7]
Structure of T_L=[7]

N=953=[953,1]
Structure of C_L=[7]
Structure of T_L=[]

N=967=[967,1]
Structure of C_L=[7]
Structure of T_L=[]

N=1009=[1009,1]
Structure of C_L=[7]
Structure of T_L=[]

N=1051=[1051,1]
Structure of C_L=[7]
Structure of T_L=[]

N=1093=[1093,1]
Structure of C_L=[7]
Structure of T_L=[]

N=1163=[1163,1]
Structure of C_L=[7]
Structure of T_L=[]

N=1289=[1289,1]
Structure of C_L=[7,7]
Structure of T_L=[7]

N=1303=[1303,1]
Structure of C_L=[7,7]
Structure of T_L=[7]

N=1373=[1373,1]
Structure of C_L=[7]
Structure of T_L=[49]

N=1429=[1429,1]
Structure of C_L=[7]
Structure of T_L=[]

N=1471=[1471,1]
Structure of C_L=[7]
Structure of T_L=[49]

N=1499=[1499,1]
Structure of C_L=[7]
Structure of T_L=[7]

N=1583=[1583,1]
Structure of C_L=[7]
Structure of T_L=[7]

N=1597=[1597,1]
Structure of C_L=[7]
Structure of T_L=[]

N=1667=[1667,1]
Structure of C_L=[7]
Structure of T_L=[49]

N=1709=[1709,1]
Structure of C_L=[7]
Structure of T_L=[]

N=1723=[1723,1]
Structure of C_L=[7,7]
Structure of T_L=[7]

N=1877=[1877,1]
Structure of C_L=[7]
Structure of T_L=[]

N=1933=[1933,1]
Structure of C_L=[7]
Structure of T_L=[]

N=2003=[2003,1]
Structure of C_L=[7]
Structure of T_L=[]

N=2017=[2017,1]
Structure of C_L=[7]
Structure of T_L=[]

N=2087=[2087,1]
Structure of C_L=[7]
Structure of T_L=[]

N=2129=[2129,1]
Structure of C_L=[7]
Structure of T_L=[]

N=2143=[2143,1]
Structure of C_L=[7]
Structure of T_L=[]

N=2213=[2213,1]
Structure of C_L=[7]
Structure of T_L=[]

N=2269=[2269,1]
Structure of C_L=[7]
Structure of T_L=[]

N=2297=[2297,1]
Structure of C_L=[7]
Structure of T_L=[]

N=2311=[2311,1]
Structure of C_L=[7]
Structure of T_L=[7]

N=2339=[2339,1]
Structure of C_L=[7,7]
Structure of T_L=[7]

N=2381=[2381,1]
Structure of C_L=[7]
Structure of T_L=[]

N=2423=[2423,1]
Structure of C_L=[7]
Structure of T_L=[7]

N=2437=[2437,1]
Structure of C_L=[7]
Structure of T_L=[]

N=2521=[2521,1]
Structure of C_L=[7]
Structure of T_L=[]

N=2549=[2549,1]
Structure of C_L=[7]
Structure of T_L=[49]

N=2591=[2591,1]
Structure of C_L=[7]
Structure of T_L=[]

N=2633=[2633,1]
Structure of C_L=[7,7]
Structure of T_L=[7]

N=2647=[2647,1]
Structure of C_L=[7]
Structure of T_L=[49]

N=2689=[2689,1]
Structure of C_L=[7]
Structure of T_L=[]

N=2731=[2731,1]
Structure of C_L=[7,7]
Structure of T_L=[7]

N=2801=[2801,1]
Structure of C_L=[7]
Structure of T_L=[7]

N=2843=[2843,1]
Structure of C_L=[7]
Structure of T_L=[49]

N=2857=[2857,1]
Structure of C_L=[7]
Structure of T_L=[7]

N=2927=[2927,1]
Structure of C_L=[7]
Structure of T_L=[7]

N=2969=[2969,1]
Structure of C_L=[7,7]
Structure of T_L=[7]

N=3011=[3011,1]
Structure of C_L=[7,7]
Structure of T_L=[7]

N=3067=[3067,1]
Structure of C_L=[7]
Structure of T_L=[7]

N=3109=[3109,1]
Structure of C_L=[7]
Structure of T_L=[]

N=3137=[3137,1]
Structure of C_L=[7]
Structure of T_L=[49]

N=3221=[3221,1]
Structure of C_L=[7]
Structure of T_L=[]

N=3319=[3319,1]
Structure of C_L=[7,7]
Structure of T_L=[7]

N=3347=[3347,1]
Structure of C_L=[7]
Structure of T_L=[7]

N=3361=[3361,1]
Structure of C_L=[7]
Structure of T_L=[]

N=3389=[3389,1]
Structure of C_L=[7,7]
Structure of T_L=[7]

N=3529=[3529,1]
Structure of C_L=[7]
Structure of T_L=[49]

N=3557=[3557,1]
Structure of C_L=[7]
Structure of T_L=[]

N=3571=[3571,1]
Structure of C_L=[7]
Structure of T_L=[]

N=3613=[3613,1]
Structure of C_L=[7,7]
Structure of T_L=[7]

N=3697=[3697,1]
Structure of C_L=[7,7]
Structure of T_L=[7]

N=3739=[3739,1]
Structure of C_L=[7,7]
Structure of T_L=[7]

N=3767=[3767,1]
Structure of C_L=[7]
Structure of T_L=[]

N=3823=[3823,1]
Structure of C_L=[7]
Structure of T_L=[49]

N=3851=[3851,1]
Structure of C_L=[7]
Structure of T_L=[]

N=3907=[3907,1]
Structure of C_L=[7]
Structure of T_L=[]

N=4019=[4019,1]
Structure of C_L=[7]
Structure of T_L=[49]

N=4159=[4159,1]
Structure of C_L=[7,7]
Structure of T_L=[7]

N=4201=[4201,1]
Structure of C_L=[7]
Structure of T_L=[]

N=4229=[4229,1]
Structure of C_L=[7,7]
Structure of T_L=[7]

N=4243=[4243,1]
Structure of C_L=[7]
Structure of T_L=[]

N=4271=[4271,1]
Structure of C_L=[7,7]
Structure of T_L=[7]

N=4327=[4327,1]
Structure of C_L=[7]
Structure of T_L=[]

N=4397=[4397,1]
Structure of C_L=[7]
Structure of T_L=[7]

N=4481=[4481,1]
Structure of C_L=[7,7]
Structure of T_L=[7]

N=4523=[4523,1]
Structure of C_L=[7,7]
Structure of T_L=[7]

N=4621=[4621,1]
Structure of C_L=[7]
Structure of T_L=[]

N=4649=[4649,1]
Structure of C_L=[7,7]
Structure of T_L=[7]

N=4663=[4663,1]
Structure of C_L=[7,7]
Structure of T_L=[7]

N=4691=[4691,1]
Structure of C_L=[7,7]
Structure of T_L=[7]

N=4733=[4733,1]
Structure of C_L=[7]
Structure of T_L=[]

N=4789=[4789,1]
Structure of C_L=[7]
Structure of T_L=[]

N=4817=[4817,1]
Structure of C_L=[7]
Structure of T_L=[]

N=4831=[4831,1]
Structure of C_L=[7]
Structure of T_L=[]

N=4943=[4943,1]
Structure of C_L=[7,7]
Structure of T_L=[7]

N=4957=[4957,1]
Structure of C_L=[7,7]
Structure of T_L=[7]

N=4999=[4999,1]
Structure of C_L=[7]
Structure of T_L=[49]

N=5153=[5153,1]
Structure of C_L=[7]
Structure of T_L=[]

N=5167=[5167,1]
Structure of C_L=[7]
Structure of T_L=[]

N=5209=[5209,1]
Structure of C_L=[7]
Structure of T_L=[]

N=5237=[5237,1]
Structure of C_L=[7]
Structure of T_L=[]

N=5279=[5279,1]
Structure of C_L=[7,7]
Structure of T_L=[7]

N=5419=[5419,1]
Structure of C_L=[7]
Structure of T_L=[]

N=5503=[5503,1]
Structure of C_L=[7]
Structure of T_L=[]

N=5531=[5531,1]
Structure of C_L=[7]
Structure of T_L=[7]

N=5573=[5573,1]
Structure of C_L=[7]
Structure of T_L=[]

N=5657=[5657,1]
Structure of C_L=[7]
Structure of T_L=[]

N=5741=[5741,1]
Structure of C_L=[7]
Structure of T_L=[]

N=5783=[5783,1]
Structure of C_L=[7]
Structure of T_L=[49]

N=5839=[5839,1]
Structure of C_L=[7]
Structure of T_L=[]

N=5867=[5867,1]
Structure of C_L=[7]
Structure of T_L=[]

N=5881=[5881,1]
Structure of C_L=[7]
Structure of T_L=[49]

N=5923=[5923,1]
Structure of C_L=[7]
Structure of T_L=[]

N=6007=[6007,1]
Structure of C_L=[7,7]
Structure of T_L=[7]

N=6091=[6091,1]
Structure of C_L=[7,7]
Structure of T_L=[7]

N=6133=[6133,1]
Structure of C_L=[7]
Structure of T_L=[]

N=6203=[6203,1]
Structure of C_L=[7,7]
Structure of T_L=[7]

N=6217=[6217,1]
Structure of C_L=[7]
Structure of T_L=[7]

N=6287=[6287,1]
Structure of C_L=[7]
Structure of T_L=[]

N=6301=[6301,1]
Structure of C_L=[7]
Structure of T_L=[]

N=6329=[6329,1]
Structure of C_L=[7]
Structure of T_L=[]

N=6343=[6343,1]
Structure of C_L=[7,7]
Structure of T_L=[7]

N=6427=[6427,1]
Structure of C_L=[7]
Structure of T_L=[7]

N=6469=[6469,1]
Structure of C_L=[7]
Structure of T_L=[49]

N=6553=[6553,1]
Structure of C_L=[7]
Structure of T_L=[]

N=6581=[6581,1]
Structure of C_L=[7]
Structure of T_L=[]

N=6637=[6637,1]
Structure of C_L=[7]
Structure of T_L=[]

N=6679=[6679,1]
Structure of C_L=[7,7]
Structure of T_L=[7]

N=6763=[6763,1]
Structure of C_L=[7]
Structure of T_L=[49]

N=6791=[6791,1]
Structure of C_L=[7,7]
Structure of T_L=[7]

N=6833=[6833,1]
Structure of C_L=[7]
Structure of T_L=[]

N=6917=[6917,1]
Structure of C_L=[7]
Structure of T_L=[7]

N=6959=[6959,1]
Structure of C_L=[7]
Structure of T_L=[49]

N=7001=[7001,1]
Structure of C_L=[7]
Structure of T_L=[]

N=7043=[7043,1]
Structure of C_L=[7,7]
Structure of T_L=[7]

N=7057=[7057,1]
Structure of C_L=[7,7]
Structure of T_L=[49]

N=7127=[7127,1]
Structure of C_L=[7]
Structure of T_L=[]

N=7211=[7211,1]
Structure of C_L=[7]
Structure of T_L=[]

N=7253=[7253,1]
Structure of C_L=[7,7]
Structure of T_L=[49]

N=7309=[7309,1]
Structure of C_L=[7,7,7]
Structure of T_L=[7,7]

N=7351=[7351,1]
Structure of C_L=[7]
Structure of T_L=[49]

N=7393=[7393,1]
Structure of C_L=[7]
Structure of T_L=[]

N=7477=[7477,1]
Structure of C_L=[7]
Structure of T_L=[]

N=7547=[7547,1]
Structure of C_L=[7,7]
Structure of T_L=[49]

N=7561=[7561,1]
Structure of C_L=[7]
Structure of T_L=[]

N=7589=[7589,1]
Structure of C_L=[7]
Structure of T_L=[]

N=7603=[7603,1]
Structure of C_L=[7]
Structure of T_L=[]

N=7673=[7673,1]
Structure of C_L=[7,7]
Structure of T_L=[7]

N=7687=[7687,1]
Structure of C_L=[7]
Structure of T_L=[7]

N=7757=[7757,1]
Structure of C_L=[7]
Structure of T_L=[]

N=7841=[7841,1]
Structure of C_L=[7]
Structure of T_L=[49]

N=7883=[7883,1]
Structure of C_L=[7,7]
Structure of T_L=[7]

N=8009=[8009,1]
Structure of C_L=[7]
Structure of T_L=[]

N=8093=[8093,1]
Structure of C_L=[7]
Structure of T_L=[]

N=8191=[8191,1]
Structure of C_L=[7]
Structure of T_L=[]

N=8219=[8219,1]
Structure of C_L=[7]
Structure of T_L=[]

N=8233=[8233,1]
Structure of C_L=[7]
Structure of T_L=[49]

N=8317=[8317,1]
Structure of C_L=[7]
Structure of T_L=[]

N=8387=[8387,1]
Structure of C_L=[7,7]
Structure of T_L=[7]

N=8429=[8429,1]
Structure of C_L=[7]
Structure of T_L=[49]

N=8443=[8443,1]
Structure of C_L=[7,7]
Structure of T_L=[7]

N=8513=[8513,1]
Structure of C_L=[7]
Structure of T_L=[]

N=8527=[8527,1]
Structure of C_L=[7,7]
Structure of T_L=[49]

N=8597=[8597,1]
Structure of C_L=[7]
Structure of T_L=[7]

N=8681=[8681,1]
Structure of C_L=[7]
Structure of T_L=[]

N=8737=[8737,1]
Structure of C_L=[7]
Structure of T_L=[7]

N=8779=[8779,1]
Structure of C_L=[7]
Structure of T_L=[7]

N=8807=[8807,1]
Structure of C_L=[7]
Structure of T_L=[]

N=8821=[8821,1]
Structure of C_L=[7]
Structure of T_L=[49]

N=8849=[8849,1]
Structure of C_L=[7]
Structure of T_L=[]

N=8863=[8863,1]
Structure of C_L=[7]
Structure of T_L=[]

N=8933=[8933,1]
Structure of C_L=[7]
Structure of T_L=[]

N=9059=[9059,1]
Structure of C_L=[7,7]
Structure of T_L=[7]

N=9157=[9157,1]
Structure of C_L=[7,7]
Structure of T_L=[7]

N=9199=[9199,1]
Structure of C_L=[7,7]
Structure of T_L=[7]

N=9227=[9227,1]
Structure of C_L=[7]
Structure of T_L=[]

N=9241=[9241,1]
Structure of C_L=[7]
Structure of T_L=[]

N=9283=[9283,1]
Structure of C_L=[7]
Structure of T_L=[]

N=9311=[9311,1]
Structure of C_L=[7]
Structure of T_L=[49]

N=9437=[9437,1]
Structure of C_L=[7]
Structure of T_L=[]

N=9479=[9479,1]
Structure of C_L=[7,7]
Structure of T_L=[7]

N=9521=[9521,1]
Structure of C_L=[7]
Structure of T_L=[]

N=9619=[9619,1]
Structure of C_L=[7]
Structure of T_L=[]

N=9661=[9661,1]
Structure of C_L=[7,7]
Structure of T_L=[7]

N=9689=[9689,1]
Structure of C_L=[7,7]
Structure of T_L=[7]

N=9787=[9787,1]
Structure of C_L=[7]
Structure of T_L=[]

N=9829=[9829,1]
Structure of C_L=[7]
Structure of T_L=[]

N=9857=[9857,1]
Structure of C_L=[7]
Structure of T_L=[]

N=9871=[9871,1]
Structure of C_L=[7,7]
Structure of T_L=[7]

N=9941=[9941,1]
Structure of C_L=[7]
Structure of T_L=[]

N=10039=[10039,1]
Structure of C_L=[7]
Structure of T_L=[]

N=10067=[10067,1]
Structure of C_L=[7,7]
Structure of T_L=[7]

N=10151=[10151,1]
Structure of C_L=[7]
Structure of T_L=[]

N=10193=[10193,1]
Structure of C_L=[7]
Structure of T_L=[49]

N=10333=[10333,1]
Structure of C_L=[7]
Structure of T_L=[]

N=10459=[10459,1]
Structure of C_L=[7,7]
Structure of T_L=[7]

N=10487=[10487,1]
Structure of C_L=[7]
Structure of T_L=[49]

N=10501=[10501,1]
Structure of C_L=[7,7]
Structure of T_L=[7]

N=10529=[10529,1]
Structure of C_L=[7]
Structure of T_L=[]

N=10613=[10613,1]
Structure of C_L=[7,7]
Structure of T_L=[7]

N=10627=[10627,1]
Structure of C_L=[7,7]
Structure of T_L=[7]

N=10711=[10711,1]
Structure of C_L=[7]
Structure of T_L=[]

N=10739=[10739,1]
Structure of C_L=[7,7]
Structure of T_L=[7]

N=10753=[10753,1]
Structure of C_L=[7]
Structure of T_L=[]

N=10781=[10781,1]
Structure of C_L=[7]
Structure of T_L=[49]

N=10837=[10837,1]
Structure of C_L=[7]
Structure of T_L=[7]

N=10949=[10949,1]
Structure of C_L=[7]
Structure of T_L=[]

N=11047=[11047,1]
Structure of C_L=[7]
Structure of T_L=[]

N=11117=[11117,1]
Structure of C_L=[7]
Structure of T_L=[7]

N=11131=[11131,1]
Structure of C_L=[7]
Structure of T_L=[7]

N=11159=[11159,1]
Structure of C_L=[7,7]
Structure of T_L=[7]

N=11173=[11173,1]
Structure of C_L=[7]
Structure of T_L=[49]

N=11243=[11243,1]
Structure of C_L=[7,7]
Structure of T_L=[7]

N=11257=[11257,1]
Structure of C_L=[7]
Structure of T_L=[]

N=11299=[11299,1]
Structure of C_L=[7,7]
Structure of T_L=[7]

N=11369=[11369,1]
Structure of C_L=[7]
Structure of T_L=[49]

N=11383=[11383,1]
Structure of C_L=[7,7]
Structure of T_L=[7]

N=11411=[11411,1]
Structure of C_L=[7]
Structure of T_L=[]

N=11467=[11467,1]
Structure of C_L=[7]
Structure of T_L=[49]

N=11551=[11551,1]
Structure of C_L=[7]
Structure of T_L=[]

N=11579=[11579,1]
Structure of C_L=[7,7]
Structure of T_L=[7]

N=11593=[11593,1]
Structure of C_L=[7]
Structure of T_L=[7]

N=11621=[11621,1]
Structure of C_L=[7]
Structure of T_L=[7]

N=11677=[11677,1]
Structure of C_L=[7]
Structure of T_L=[]

N=11719=[11719,1]
Structure of C_L=[7,7]
Structure of T_L=[7]

N=11789=[11789,1]
Structure of C_L=[7,7]
Structure of T_L=[7]

N=11831=[11831,1]
Structure of C_L=[7]
Structure of T_L=[]

N=11887=[11887,1]
Structure of C_L=[7]
Structure of T_L=[]

N=11971=[11971,1]
Structure of C_L=[7]
Structure of T_L=[]

N=12041=[12041,1]
Structure of C_L=[7]
Structure of T_L=[]

N=12097=[12097,1]
Structure of C_L=[7]
Structure of T_L=[]

N=12251=[12251,1]
Structure of C_L=[7]
Structure of T_L=[49]

N=12377=[12377,1]
Structure of C_L=[7,7,7]
Structure of T_L=[7,7]

N=12391=[12391,1]
Structure of C_L=[7]
Structure of T_L=[]

N=12433=[12433,1]
Structure of C_L=[7]
Structure of T_L=[]

N=12503=[12503,1]
Structure of C_L=[7]
Structure of T_L=[]

N=12517=[12517,1]
Structure of C_L=[7,7]
Structure of T_L=[7]

N=12601=[12601,1]
Structure of C_L=[7]
Structure of T_L=[]

N=12671=[12671,1]
Structure of C_L=[7,7]
Structure of T_L=[7]

N=12713=[12713,1]
Structure of C_L=[7,7]
Structure of T_L=[7]

N=12853=[12853,1]
Structure of C_L=[7]
Structure of T_L=[]

N=12923=[12923,1]
Structure of C_L=[7]
Structure of T_L=[]

N=12979=[12979,1]
Structure of C_L=[7]
Structure of T_L=[]

N=13007=[13007,1]
Structure of C_L=[7,7]
Structure of T_L=[7]

N=13049=[13049,1]
Structure of C_L=[7]
Structure of T_L=[]

N=13063=[13063,1]
Structure of C_L=[7,7]
Structure of T_L=[7]

N=13147=[13147,1]
Structure of C_L=[7,7]
Structure of T_L=[7]

N=13217=[13217,1]
Structure of C_L=[7,7]
Structure of T_L=[7]

N=13259=[13259,1]
Structure of C_L=[7]
Structure of T_L=[]

N=13399=[13399,1]
Structure of C_L=[7]
Structure of T_L=[]

N=13441=[13441,1]
Structure of C_L=[7,7,7]
Structure of T_L=[7,7]

N=13469=[13469,1]
Structure of C_L=[7,7]
Structure of T_L=[7]

N=13553=[13553,1]
Structure of C_L=[7]
Structure of T_L=[7]

N=13567=[13567,1]
Structure of C_L=[7]
Structure of T_L=[]

N=13679=[13679,1]
Structure of C_L=[7,7]
Structure of T_L=[7]

N=13693=[13693,1]
Structure of C_L=[7]
Structure of T_L=[]

N=13721=[13721,1]
Structure of C_L=[7]
Structure of T_L=[49]

N=13763=[13763,1]
Structure of C_L=[7,7]
Structure of T_L=[7]

N=13903=[13903,1]
Structure of C_L=[7]
Structure of T_L=[]

N=13931=[13931,1]
Structure of C_L=[7]
Structure of T_L=[]

N=14029=[14029,1]
Structure of C_L=[7]
Structure of T_L=[]

N=14057=[14057,1]
Structure of C_L=[7]
Structure of T_L=[]

N=14071=[14071,1]
Structure of C_L=[7]
Structure of T_L=[]

N=14197=[14197,1]
Structure of C_L=[7,7,7,7]
Structure of T_L=[7,7,7]

N=14281=[14281,1]
Structure of C_L=[7]
Structure of T_L=[]

N=14323=[14323,1]
Structure of C_L=[7,7]
Structure of T_L=[7]

N=14407=[14407,1]
Structure of C_L=[7,7]
Structure of T_L=[49]

N=14449=[14449,1]
Structure of C_L=[7,7]
Structure of T_L=[7]

N=14519=[14519,1]
Structure of C_L=[7,7]
Structure of T_L=[7]

N=14533=[14533,1]
Structure of C_L=[7,7]
Structure of T_L=[7]

N=14561=[14561,1]
Structure of C_L=[7,7]
Structure of T_L=[7]
\end{verbatim}
\normalsize

\onecolumn
\subsection{Table of $\Cl_L$, ${\mathcal T}_L$, $p=11$ 
with arbitrary $N$}
${}$
\scriptsize
\begin{verbatim}
N=11=[11,1]
Structure of C_L=[]
Structure of T_L=[]

N=23=[23,1]
Structure of C_L=[11]
Structure of T_L=[]

N=67=[67,1]
Structure of C_L=[11,11]
Structure of T_L=[11]

N=89=[89,1]
Structure of C_L=[11]
Structure of T_L=[]

N=199=[199,1]
Structure of C_L=[11,11]
Structure of T_L=[11]

N=253=[11,1;23,1]
Structure of C_L=[11]
Structure of T_L=[]

N=331=[331,1]
Structure of C_L=[11]
Structure of T_L=[]

N=353=[353,1]
Structure of C_L=[11]
Structure of T_L=[]

N=397=[397,1]
Structure of C_L=[11]
Structure of T_L=[]

N=419=[419,1]
Structure of C_L=[11]
Structure of T_L=[]

N=463=[463,1]
Structure of C_L=[11]
Structure of T_L=[]

N=529=[23,2]
Structure of C_L=[11]
Structure of T_L=[]

N=617=[617,1]
Structure of C_L=[11]
Structure of T_L=[]

N=661=[661,1]
Structure of C_L=[11,11]
Structure of T_L=[11]

N=683=[683,1]
Structure of C_L=[11,11]
Structure of T_L=[11]

N=727=[727,1]
Structure of C_L=[11]
Structure of T_L=[121]

N=737=[11,1;67,1]
Structure of C_L=[11]
Structure of T_L=[]
\end{verbatim}
\normalsize

\twocolumn\scriptsize
\begin{verbatim}
N=859=[859,1]
Structure of C_L=[11]
Structure of T_L=[]

N=881=[881,1]
Structure of C_L=[11]
Structure of T_L=[]

N=947=[947,1]
Structure of C_L=[11]
Structure of T_L=[]

N=979=[11,1;89,1]
Structure of C_L=[11]
Structure of T_L=[]

N=991=[991,1]
Structure of C_L=[11]
Structure of T_L=[]

N=1013=[1013,1]
Structure of C_L=[11]
Structure of T_L=[11]

N=1123=[1123,1]
Structure of C_L=[11]
Structure of T_L=[]

N=1277=[1277,1]
Structure of C_L=[11]
Structure of T_L=[]

N=1321=[1321,1]
Structure of C_L=[11,11]
Structure of T_L=[11]

N=1409=[1409,1]
Structure of C_L=[11]
Structure of T_L=[]

N=1453=[1453,1]
Structure of C_L=[11]
Structure of T_L=[121]

N=1541=[23,1;67,1]
Structure of C_L=[11,11]
Structure of T_L=[11]

N=1607=[1607,1]
Structure of C_L=[11]
Structure of T_L=[]

N=1783=[1783,1]
Structure of C_L=[11]
Structure of T_L=[]

N=1871=[1871,1]
Structure of C_L=[11]
Structure of T_L=[]

N=2003=[2003,1]
Structure of C_L=[11,11]
Structure of T_L=[11]

N=2047=[23,1;89,1]
Structure of C_L=[11,11,11]
Structure of T_L=[11,11]

N=2069=[2069,1]
Structure of C_L=[11]
Structure of T_L=[11]

N=2113=[2113,1]
Structure of C_L=[11]
Structure of T_L=[]

N=2179=[2179,1]
Structure of C_L=[11]
Structure of T_L=[121]

N=2189=[11,1;199,1]
Structure of C_L=[11]
Structure of T_L=[]

N=2267=[2267,1]
Structure of C_L=[11]
Structure of T_L=[]

N=2311=[2311,1]
Structure of C_L=[11]
Structure of T_L=[]

N=2333=[2333,1]
Structure of C_L=[11]
Structure of T_L=[]

N=2377=[2377,1]
Structure of C_L=[11,11]
Structure of T_L=[11]

N=2399=[2399,1]
Structure of C_L=[11]
Structure of T_L=[]

N=2531=[2531,1]
Structure of C_L=[11]
Structure of T_L=[11]
\end{verbatim}
\normalsize

\onecolumn
\subsection{Table of $\Cl_L$, ${\mathcal T}_L$, $p=11$, 
$N$ prime, $N \equiv 1 \pmod p$}
${}$
\scriptsize
\begin{verbatim}
N=2663=[2663,1]
Structure of C_L=[11]
Structure of T_L=[121]

N=2707=[2707,1]
Structure of C_L=[11,11]
Structure of T_L=[11]

N=2729=[2729,1]
Structure of C_L=[11]
Structure of T_L=[]

N=2861=[2861,1]
Structure of C_L=[11]
Structure of T_L=[]

N=2927=[2927,1]
Structure of C_L=[11]
Structure of T_L=[]

N=2971=[2971,1]
Structure of C_L=[11]
Structure of T_L=[]

N=3037=[3037,1]
Structure of C_L=[11]
Structure of T_L=[]

N=3169=[3169,1]
Structure of C_L=[11]
Structure of T_L=[]

N=3191=[3191,1]
Structure of C_L=[11,11]
Structure of T_L=[11]

N=3257=[3257,1]
Structure of C_L=[11]
Structure of T_L=[]

N=3301=[3301,1]
Structure of C_L=[11]
Structure of T_L=[]

N=3323=[3323,1]
Structure of C_L=[11]
Structure of T_L=[11]

N=3389=[3389,1]
Structure of C_L=[11]
Structure of T_L=[121]

N=3433=[3433,1]
Structure of C_L=[11]
Structure of T_L=[]

N=3499=[3499,1]
Structure of C_L=[11]
Structure of T_L=[]

N=3631=[3631,1]
Structure of C_L=[11,11]
Structure of T_L=[121]

N=3697=[3697,1]
Structure of C_L=[11]
Structure of T_L=[]

N=3719=[3719,1]
Structure of C_L=[11]
Structure of T_L=[11]

N=3851=[3851,1]
Structure of C_L=[11]
Structure of T_L=[]

N=3917=[3917,1]
Structure of C_L=[11]
Structure of T_L=[]

N=4027=[4027,1]
Structure of C_L=[11,11]
Structure of T_L=[11]

N=4049=[4049,1]
Structure of C_L=[11]
Structure of T_L=[]

N=4093=[4093,1]
Structure of C_L=[11]
Structure of T_L=[]

N=4159=[4159,1]
Structure of C_L=[11]
Structure of T_L=[]

N=4357=[4357,1]
Structure of C_L=[11]
Structure of T_L=[121]

N=4423=[4423,1]
Structure of C_L=[11]
Structure of T_L=[]

N=4621=[4621,1]
Structure of C_L=[11]
Structure of T_L=[]

N=4643=[4643,1]
Structure of C_L=[11]
Structure of T_L=[]

N=4951=[4951,1]
Structure of C_L=[11]
Structure of T_L=[]

N=4973=[4973,1]
Structure of C_L=[11]
Structure of T_L=[]

N=5039=[5039,1]
Structure of C_L=[11]
Structure of T_L=[]
\end{verbatim}
\normalsize

\onecolumn
\subsection{Table of $\Cl_L$, ${\mathcal T}_L$, $p=13$, 
$N$ prime, $N \equiv 1 \pmod p$}
${}$
\scriptsize
\begin{verbatim}
N=53=[53,1]
Structure of C_L=[13]
Structure of T_L=[]

N=79=[79,1]
Structure of C_L=[13]
Structure of T_L=[]

N=131=[131,1]
Structure of C_L=[13]
Structure of T_L=[]

N=157=[157,1]
Structure of C_L=[13]
Structure of T_L=[13]

N=313=[313,1]
Structure of C_L=[13]
Structure of T_L=[]

N=443=[443,1]
Structure of C_L=[13]
Structure of T_L=[13]

N=521=[521,1]
Structure of C_L=[13]
Structure of T_L=[]

N=547=[547,1]
Structure of C_L=[13]
Structure of T_L=[13]

N=599=[599,1]
Structure of C_L=[13]
Structure of T_L=[]

N=677=[677,1]
Structure of C_L=[13]
Structure of T_L=[169]

N=859=[859,1]
Structure of C_L=[13]
Structure of T_L=[13]

N=911=[911,1]
Structure of C_L=[13]
Structure of T_L=[]

N=937=[937,1]
Structure of C_L=[13]
Structure of T_L=[13]

N=1093=[1093,1]
Structure of C_L=[13]
Structure of T_L=[]

N=1171=[1171,1]
Structure of C_L=[13]
Structure of T_L=[]

N=1223=[1223,1]
Structure of C_L=[13]
Structure of T_L=[]

N=1249=[1249,1]
Structure of C_L=[13]
Structure of T_L=[]

N=1301=[1301,1]
Structure of C_L=[13]
Structure of T_L=[]

N=1327=[1327,1]
Structure of C_L=[13]
Structure of T_L=[]

N=1483=[1483,1]
Structure of C_L=[13,13]
Structure of T_L=[13]

N=1613=[1613,1]
Structure of C_L=[13]
Structure of T_L=[]

N=1847=[1847,1]
Structure of C_L=[13]
Structure of T_L=[]

N=1873=[1873,1]
Structure of C_L=[13]
Structure of T_L=[]

N=1951=[1951,1]
Structure of C_L=[13]
Structure of T_L=[13]

N=2003=[2003,1]
Structure of C_L=[13]
Structure of T_L=[]

N=2029=[2029,1]
Structure of C_L=[13]
Structure of T_L=[169]

N=2081=[2081,1]
Structure of C_L=[13]
Structure of T_L=[]

N=2237=[2237,1]
Structure of C_L=[13]
Structure of T_L=[]

N=2341=[2341,1]
Structure of C_L=[13]
Structure of T_L=[]

N=2393=[2393,1]
Structure of C_L=[13]
Structure of T_L=[]

N=2549=[2549,1]
Structure of C_L=[13]
Structure of T_L=[]

N=2731=[2731,1]
Structure of C_L=[13]
Structure of T_L=[]

N=2861=[2861,1]
Structure of C_L=[13]
Structure of T_L=[13]

N=2887=[2887,1]
Structure of C_L=[13,13]
Structure of T_L=[13]
\end{verbatim}
\normalsize

\onecolumn
\section{Table of the whole class groups, $p \geq 3$} \label{C}
The idea of this table is due to some comments by Christian 
Maire to me, from a lecture given by Ren\'e Schoof in China (Harbin 2019) 
about the properties of the structure of these class groups mentioned in 
\S\,\ref{HWS}.

\smallskip
The following program is a simple adaptation of the main previous one.
It gives the structures of the $q$-class groups for any prime $q$ dividing
the class numbers. One may vary $p$ at will. 

\smallskip
For $p=5$, at the end of the 
table, we give also the invariants of the intermediate field $\Q(\sqrt 5)$,
which allows the computation of the $\omega^2$-components.

\medskip
\footnotesize
\begin{verbatim}
{p=3;PK=polcyclo(p);for(N=2,10^4,F=factor(N);
F1=component(F,1);F2=component(F,2);d=matsize(F1)[1];
T=1;for(j=1,d,e=F2[j];if(e>=p,T=0));if(T==0,next);
T=1;for(k=2,p-1,NN=1;for(j=1,d,ek=lift(Mod(k*F2[j],p));
NN=NN*F1[j]^ek);if(NN<N,T=0;break));if(T==0,next);
print();
if(Mod(N,p)!=0 & Mod(N^(p-1),p^2)==1,
print("N=",N,"=",factor(N)," ",p," unramified"));
if(Mod(N,p)==0 || Mod(N^(p-1),p^2)!=1,
print("N=",N,"=",factor(N)," ",p," ramified"));
PL=x^p-N;L=bnfinit(PL,1);
ClassL=L.cyc;\\ class number of L
Lno=L.no;Lcyc=L.cyc;\\ global structure of ClassL
DL=matsize(Lcyc)[2];FL=factor(Lno);FL1=component(FL,1);
dL=matsize(FL1)[1];
for(j=1,dL,q=FL1[j];HL=List;
for(k=1,DL,c=Lcyc[k];
if(Mod(c,q)==0,listinsert(HL,q^valuation(c,q),1)));
  \\ q-structure of ClassL
print("q=",q," ClassqL=",HL));
PM=polcompositum(PL,PK)[1];M=bnfinit(PM,1);
ClassM=M.cyc;\\ class number of M
Mno=M.no;Mcyc=M.cyc;\\ global structure of ClassM
DM=matsize(Mcyc)[2];FM=factor(Mno);FM1=component(FM,1);
dM=matsize(FM1)[1];
for(j=1,dM,q=FM1[j];HM=List;
for(k=1,DM,c=Mcyc[k];
if(Mod(c,q)==0,listinsert(HM,q^valuation(c,q),1)));
  \\ q-structure of ClassM
print("q=",q," ClassqM= ",HM)))}
\end{verbatim}
\normalsize

\twocolumn
\footnotesize

\normalsize
\textwidth=30pc

\newpage
\onecolumn


\begin{thebibliography}{xx}

\bibitem{AAIMT}S. Aouissi, A. Azizi, M.C. Ismaili, D.C. Mayer, M. Talbi,
{\it Structure of relative genus fields of cubic Kummer extensions}.
\url{https://arxiv.org/abs/1808.04678}

\bibitem{ATIA1} S. Aouissi, M. Talbi, M.C. Ismaili, A. Azizi,
{\it On a Conjecture of Lemmermeyer}.\par
\url{https://arxiv.org/abs/1810.07172}

\bibitem{ATIA2} S. Aouissi, M. Talbi, M.C. Ismaili, A. Azizi,
{\it Fields $\Q(\sqrt[3]{d},\zeta_3)$ whose $3$-class group is of type $(9,3)$},
\url{https://arxiv.org/abs/1805.04963}

\bibitem{BJ} K. Belabas,  J-F. Jaulent, {\it The logarithmic class group 
package in PARI/GP}, Pub. Math. Besan\c con 
(Th\'eorie des Nombres) (2016), 5--18.\par
\url{http://pmb.univ-fcomte.fr/2016/pmb_2016.pdf}

\bibitem{BP} F. Bertrandias, J-J. Payan, {\it $\Gamma$-extensions et invariants
cyclotomiques}, Ann. Sci. Ec. Norm. Sup. 4e s\'erie, {\bf 5}(4) (1972), 517--548.
\url{https://doi.org/10.24033/asens.1236}

\bibitem{Calegari-Emerton} F. Calegari, M. Emerton,
{\it  On the ramification of Hecke algebras at Eisenstein primes},
Invent. Math., {\bf 160}(1) (2005), 97--144.\url{https://arxiv.org/abs/math/0311368}

\bibitem{Gerth} F. Gerth III, {\it On 3-class groups of certain pure cubic fields},
Bull. Austral. Math. Soc., {\bf 72}(3) (20055), 471--476.
\url{https://doi.org/10.1017/S0004972700035292}

\bibitem{Gras0} G. Gras, {\it Sur le $3$-rang des corps cubiques non
galoisiens}, Pub. Math. Besan\c con (Th\'eorie des Nombres),
Ann\'ees 1974/1975.\url{http://pmb.univ-fcomte.fr/1975/GGras.pdf}

\bibitem{Gras1} G. Gras, {\it Class Field Theory: from theory to practice}, corr. 
2nd ed., Springer Monographs in Mathematics, Springer (2005), xiii+507 pages.\par
\url{https://doi.org/10.1007/978-3-662-11323-3}

\bibitem{Gras2}  G. Gras, {\it The $p$-adic Kummer-Leopoldt Constant:
Normalized $p$-adic Regulator}, Int. J. Number Theory {\bf 14}(2) (2018), 329--337.  
\url{https://doi.org/10.1142/S1793042118500203}

\bibitem{Gras10}  G. Gras, {\it Heuristics and conjectures in the
direction of a $p$-adic Brauer--Siegel theorem}, Math. Comp. {\bf 88}(318) (2019),
1929--1965. \url{https://doi.org/10.1090/mcom/3395}

\bibitem{Gras3} G. Gras, {\it Practice of incomplete p-ramification over 
a number field -- History of abelian p-ramification}.
\url{https://arxiv.org/pdf/1904.10707}

\bibitem{Gras4} G. Gras, {\it Th\'eor\`emes de r\'eflexion},
J. Th\'eor. Nombres Bordeaux {\bf 10}(2) (1998), 399--499.
\url{http://www.numdam.org/item/JTNB_1998__10_2_399_0/}

\bibitem{Gras5} G. Gras,  {\it On $p$-rationality of number fields. 
Applications--PARI/GP programs}, Pub. Math. Besan\c con 
(Th\'eorie des Nombres), Ann\'ees 2017/2018.\par
\url{https://arxiv.org/pdf/1709.06388}

\bibitem{Gras6}  G. Gras, {\it  Groupe de Galois de la $p$-extension 
ab\'elienne $p$-ramifi\'ee maximale d'un corps de nombres}, J. reine angew. 
Math. {\bf 333} (1982), 86--132. \par \url{https://eudml.org/doc/152440}
\url{https://www.researchgate.net/publication/243110955} 

\bibitem{Gras7}  G. Gras, {\it  Logarithme $p$-adique et groupes 
de Galois}, J. reine angew. Math. {\bf 343} (1983), 64--80.        
\url{https://doi.org/10.1515/crll.1983.343.64}

\bibitem{Gras8}  G. Gras, {\it Invariant generalized ideal 
classes--Structure theorems for $p$-class groups in $p$-extensions},
Proc. Indian Acad. Sci. (Math. Sci.) {\bf 127}(1) ( 2017), 1--34. \par
\url{https://www.ias.ac.in/article/fulltext/pmsc/127/01/0001-0034}

\bibitem{GJ}  G. Gras,  J-F. Jaulent,  {\it Sur les corps de nombres r\'eguliers}, 
Math. Z. {\bf 202}(3) (1989), 343--365.
\url{https://eudml.org/doc/174095}

\bibitem{HW} D. Hubbard , L.C. Washington, {\it Iwasawa invariants 
of some non-cyclotomic $\Z_p$-extensions}, J. Number Theory, {\bf 188},
(2018), 18-47.

\bibitem{Iimura} K. Iimura,
{\it On the $l$-rank of ideal class groups of certain number fields},
Acta Arith., {\bf 47}(2) (1986), 153--166.
\url{http://matwbn.icm.edu.pl/ksiazki/aa/aa47/aa4724.pdf}

\bibitem{Jaulent} J-F. Jaulent, {\it Unit\'es et classes dans les extensions 
m\'etab\'eliennes de degr\'e $nl^s$ sur un corps de nombres alg\'ebriques},
Ann. Inst. Fourier (Grenoble), {\bf 31}(1) (1981 ), ix--x, 39--62.
\url{http://archive.numdam.org/article/AIF_1981__31_1_39_0.pdf}

\bibitem{Jaulent0} J-F. Jaulent, {\it Th\'eorie $\ell$-adique globale du  corps de 
classes}, J. Th\'eorie des Nombres de Bordeaux, {\bf 10}(2) (1998), 355--397.\par
\url{http://www.numdam.org/article/JTNB_1998__10_2_355_0.pdf}

\bibitem{Jalog} J-F. Jaulent, {\it Classes logarithmiques des corps 
de nombres}, J. Th\'eorie des Nombres de Bordeaux {\bf 6} (1994), 301--325. \par
\url{http://www.numdam.org/article/JTNB_1994__6_2_301_0.pdf}

\bibitem{JN} J-F. Jaulent, T. Nguyen Quang Do, {\it Corps $p$-rationnels, 
corps $p$-r\'eguliers et ramification restreinte}, J. Th\'eor. Nombres 
Bordeaux {\bf 5} (1993), 343--363. 
\url{http://www.numdam.org/article/JTNB_1993__5_2_343_0.pdf}

\bibitem{Kobayashi} S. Kobayashi,  {\it Complete determination of the $3$-rank
in pure cubic fields}, J. Math. Soc. Japan {\bf 29}(2) (1977), 373--384.\par
\url{https://projecteuclid.org/download/pdf_1/euclid.jmsj/1240433454}

\bibitem{Lecouturier} E. Lecouturier, {\it On the Galois structure of the class group
of certain Kummer extensions}, J. London Math. Soc. (2) {\bf 98} (2018), 35--58.
\url{https://doi.org/10.1112/jlms.12123}

\bibitem{Movahhedi} A. Movahhedi, {\it Sur les $p$-extensions 
des corps $p$-rationnels}, Math. Nachr. {\bf 149} (1990), 163--176.
\url{http://onlinelibrary.wiley.com/doi/10.1002/mana.19901490113/}

\bibitem{MN} A. Movahhedi, T. Nguyen Quang Do, {\it Sur l'arithm\'etique 
des corps de nombres $p$-rationnels}, S\'eminaire de Th\'eorie des Nombres, 
Paris 1987--88, Progress in Math. {\bf 81} (1990), 155--200.
\url{https://doi.org/10.1007/978-1-4612-3460-9_9}

\bibitem{MR} C. Maire, M. Rougnant, {\it A note on $p$-rational 
fields and the $abc$-conjecture} (2019), Submitted to Proceedings of the AMS.
\url{https://arxiv.org/abs/1903.11271} 

\bibitem{pari} The PARI Group,  PARI/GP, version \texttt{2.9.0}, 
Universit\'e de Bordeaux (2016).\par
\url{http://pari.math.u-bordeaux.fr/}

\bibitem{Schaefer-Stubley} K. Schaefer, E. Stubley,
{\it Class groups of Kummer extensions via cup products in Galois cohomology}.
\url{https://arxiv.org/abs/1806.00517}

\bibitem{Schoof} R. Schoof, {\it Kummer extensions of cyclotomic fields}, 
Leiden, March 2018.

\bibitem{Washington} L.C. Washington, {\it Introduction to cyclotomic fields}, 
Graduate Texts in Mathematics, {\bf 83}, Springer-Verlag, 
New York, 1997, {\rm xiv}+487 pp.

\end{thebibliography}
\end{document}